\documentclass[12pt]{article}
\usepackage{latexsym,amsmath,amssymb}

\newif\ifdviwin

\usepackage[latin1]{inputenc}
\usepackage[english]{babel}
\usepackage{indentfirst}
\usepackage[mathscr]{eucal}
\usepackage{amssymb,amsmath,amsfonts}
\usepackage{fancybox,fancyhdr}
\usepackage{graphicx}

\newif\ifdviwin

\dviwintrue

\def\cA{\mathcal{A}}

\def\cS{\mathcal{S}}
\def\cC{\mathcal{C}}
\def\cR{\mathcal{R}}
\def\cB{\mathcal{B}}
\def\cT{\mathcal{T}}

\def\cU{\mathcal{U}}
\def\cG{\mathcal{G}}

\def\cT{\mathcal{T}}

\let\8=\infty \let\0=\emptyset  
\let\hat=\widehat

\let\tilde=\widetilde

\let\landa=\lambda
\let\alfa=\alpha

\let\parc=\partial

\def\landa{\lambda}
\def\Landa{\Lambda}

\def\flecha{\rightarrow}
\def\esiz{\langle}
\def\esde{\rangle}

\def\cte.{\mathop{\rm cte.}\nolimits}

\def\det{\mathop{\rm det}\nolimits}

\def\Re{\mathop{\rm Re }\nolimits} 

\def\sen{\mathop{\rm sen }\nolimits}

\def\N{\mathbb{N}}

\def\Z{\mathbb{Z}}
\def\Q{\mathbb{Q}}
\def\R{\mathbb{R}}
\def\T{\mathbb{T}}
\def\C{\mathbb{C}}

\def\S{\mathbb{S}}
\def\X{\mathfrak{X}}
 \newtheorem{defi}{Definition}
 \newtheorem{teo}[defi]{Theorem}
 \newtheorem{pro}[defi]{Proposition}
 \newtheorem{cor}[defi]{Corollary}
 \newtheorem{lem}[defi]{Lemma}
 
 \newtheorem{remark}[defi]{Remark}

 \newenvironment{proof}{\rm \trivlist \item[\hskip \labelsep{\it
      Proof}:]}{\par\nopagebreak \hfill $\Box$ \endtrivlist}

 \newenvironment{proof1}{\rm \trivlist \item[\hskip \labelsep{\it
      Proof of Theorem \ref{Boris}}:]}{\par\nopagebreak \hfill $\Box$ \endtrivlist}

 \newenvironment{proof2}{\rm \trivlist \item[\hskip \labelsep{\it
      Proof of Theorem \ref{truco}}:]}{\par\nopagebreak \hfill $\Box$ \endtrivlist}

\begin{document}

\vspace{1cm}

\mbox{}\vspace{0.1cm}

\begin{center}
\rule{14cm}{1.5pt}\vspace{0.5cm}

{\Large \bf Isometric immersions of $\R^2$ into $\R^4$ and
\\\vspace{0.3cm} perturbation of Hopf tori}\\ \vspace{0.5cm} {\large José A. Gálvez$\mbox{}^a$
and Pablo Mira$\mbox{}^b$}\\ \vspace{0.3cm} \rule{14cm}{1.5pt}
\end{center}
 \vspace{1cm}
$\mbox{}^a$ Departamento de Geometría y Topología, Universidad de Granada,
E-18071 Granada, Spain. \\ e-mail: jagalvez@ugr.es \vspace{0.2cm}

\noindent $\mbox{}^b$ Departamento de Matemática Aplicada y Estadística,
Universidad Politécnica de Cartagena, E-30203 Cartagena, Murcia, Spain. \\
e-mail: pablo.mira@upct.es \vspace{0.3cm}

 \begin{abstract}
We produce a new general family of flat tori in $\R^4$, the first one
since Bianchi's classical works in the 19th century. To construct these
flat tori, obtained via small perturbation of certain Hopf tori in $\S^3$,
we first present a global description of all isometric immersions of
$\R^2$ into $\R^4$ with flat normal bundle.
 \end{abstract}
\section{Introduction}
Isometric immersions of $\R^2$ into $\R^4$ appear as a critical
situation in one of the central problems of the differential geometry
of submanifolds, namely, the problem of classifying isometric
immersions between space forms. Let $Q^n (c),\bar{Q}^{n+p} (c)$ be
two space forms of dimensions $n,n+p$ respectively, and with the same
constant curvature $c$. If $c\geq 0$ and $p<n$ there exist very few
isometric immersions $f:Q^n\flecha \bar{Q}^{n+p}$, which can be
thought of as {\it trivial}. However, if
$c\geq 0$ and $p\geq n$, or if $c<0$, the situation changes completely and there exists a wide variety of
isometric immersions from $Q^n$ into $\bar{Q}^{n+p}$. This leaves the
case
$c=0$ and $n=p=2$, that is, isometric immersions of $\R^2$ into $\R^4$, as the limit
case between both situations.

There is another fact that motivates the study of isometric
immersions of $\R^2$ into $\R^4$. Recall that a complete flat
surface must be homeomorphic to a plane, a cylinder, a torus, a
Möbius strip or a Klein bottle. Among these five topological
cases, only the plane and the cylinder can be realized as complete
flat surfaces in $\R^3$. On the other hand, all of them occur when
one considers flat surfaces in $\R^4$. This indicates that $\R^4$
is the most natural ambient space in where to consider a flat
surface. One further motivation in connection with this point is
given by Tompkins theorem, which states that there are no compact
flat $n$-manifolds isometrically immersed in $\R^{2n-1}$. Again,
the existence of flat tori in $\R^4$ makes isometric immersions of
$\R^2$ into $\R^4$ a critical situation.

Motivated by all of this, the present paper studies isometric
immersions of $\R^2$ into $\R^4$ with \emph{flat normal bundle}, i.e.
$R^{\perp}\equiv 0$. This geometric setting contains a large number of
examples, and has been widely studied (see
\cite{Bia,Bor,CaDa,DaTo,DaTo2,Eno2,HKP,Kit1,Kit2,Spi,Ten,Wei1,Wei2,Wei3})).
We note that the flat normal bundle condition is the most natural
hypothesis when dealing with isometric immersions between space forms
with large codimension, as it is our case.

In this work we give a global description of all simply connected
flat surfaces with flat normal bundle in $\R^4$, and we use it as
main application to study flat tori in $\R^4$ with $R^{\perp}\equiv
0$ and regular Gauss map.

Up to now, and except for some isolated examples with
$R^{\perp}\neq 0$, only two families of flat tori in $\R^4$
have been found, and both were already known in the 19th century. One
of these families is made up by the \emph{product tori}
$\gamma_1\times\gamma_2 \subset\R^2\times \R^2\equiv\R^4$, where
$\gamma_1,\gamma_2$ are regular closed curves in $\R^2$. The
other family consists of flat tori in the $3$-sphere $\S^3$. In
\cite{Bia} Bianchi showed that the multiplication in the quaternionic
model of $\S^3$ of two closed curves in
$\S^3$ satisfying certain conditions gives rise to a flat torus in $\S^3$ for which
these curves are asymptotic curves. Much later, Kitagawa proved in
\cite{Kit1} that the asymptotic curves of flat tori in $\S^3$ are
periodic, and thus the classical Bianchi method describes all flat
tori in $\S^3$. In addition, Weiner \cite{Wei1} classified all flat
tori in $\S^3$ by means of their Gauss maps, thus solving an open
problem posed by S.T. Yau \cite{Yau}.

In this way, no new flat torus in $\R^4$ with $R^{\perp}\equiv 0$ has
been discovered since the 19th century, even though this family has
been deeply studied. Besides, Weiner showed in \cite{Wei1} that there
exist no flat Klein bottles in $\R^4$ with $R^{\perp}\equiv 0$ and
regular Gauss map. These facts created the impression, explicitly
stated as a conjecture in \cite{Bor}, that there do not exist flat
tori in $\R^4$ with $R^{\perp}\equiv 0$ apart from the previous two
families.

The main result of the present work makes use of the above mentioned
global des\-cription of flat surfaces with
$R^{\perp}\equiv 0$ in $\R^4$ to construct infinitely many new
flat tori in $\R^4$ with flat normal bundle and regular Gauss map,
that do not lie in any of these two families. This produces the
first general family of flat tori in $\R^4$ that comes out since
the 19th century.

We have organized this paper as follows. In Section $2$ we study
the Gauss map of flat surfaces $f:\Sigma\flecha\R^4$ with
$R^{\perp}\equiv 0$. For that purpose, we introduce the notion of
\emph{flat map}, which is essentially a map $F:\Sigma\flecha\S^3$
such that at its regular points it is a flat surface in $\S^3$.
Through this concept it is shown that if the Gauss map $\cG$ of
$f$ is regular, the surface $\Sigma$ inherits via $\cG$ a
canonical Lorentz surface structure. Moreover, there exists an
orthonormal basis $\{N,\hat{N}\}$ of $\X^{\perp}(f)$ such that
both $N,\hat{N}$ are flat maps whose asymptotic directions at
their regular points coincide with the null directions of the
Lorentz surface $\Sigma$. This fact is used to endow $\Sigma$ with
a global coordinate immersion that generalizes the well known
existence of Tschebysheff coordinates on any complete flat surface
in $\S^3$. As applications, we show that global Tschebysheff
coordinates exist on any simply connected \emph{analytic} flat
surface in $\S^3$, and that every analytic flat surface in $\S^3$
is orientable, regardless of its completeness. In contrast, we
produce $C^{\8}$ flat Möbius strips in $\S^3$ via the Hopf
fibration.

In Section $3$ we obtain a global representation for all simply connected
flat surfaces in $\R^4$ with flat normal bundle and regular Gauss map, in
terms of flat maps. This representation depends on the integration of a
hyperbolic linear differential system, of which we present some partial
integrations, by using both analytic and geometric methods. It is
important to observe that, for the local situation, there already exists a
representation formula due to do Carmo and Dajczer \cite{CaDa}. In there,
the authors describe (except for \emph{compositions}) all isometric
immersions of $U\subset \R^2$ into $\R^4$, with $U$ sufficiently small.

Finally, in Section $4$ we construct a large family of new flat tori in
$\R^4$, all of them with flat normal bundle and regular Gauss map. They
are not contained in any affine $3$-sphere, and cannot be expressed as the
product of two curves. To do so, we first apply the representation formula
of Section $3$ to obtain a procedure of unfolding a Hopf torus in $\S^3$,
so that one gets a flat surface in $\R^4$ (possibly with singular points)
with flat normal bundle that does not lie in any affine $3$-sphere of
$\R^4$. Then we show that for a certain family of Hopf tori this procedure
generates flat tori in $\R^4$. At last, we collapse the flat torus we have
just constructed into the Hopf torus we started with, and we get rid of
singular points. This finishes the construction. We also show in this
Section how this ideas can be used to construct new complete flat
cylinders in $\R^4$.

While this paper was being prepared, the authors knew about an interesting
paper by Weiner \cite{Wei3}, in where the author reports on his advances
in the problem of describing all isometric immersions of $\R^2$ into
$\R^4$. As a matter of fact, he mainly focuses on the $R^{\perp}\equiv 0$
case, which is the situation studied in the present work, but our approach
has important differences with respect to Weiner's one. Indeed, our
description recovers the proper immersion $f:\Sigma\flecha\R^4$, and not
just its differential $df$, and this is crucial in the construction of new
flat tori in $\R^4$ that we achieve here.

\section{Study of the Gauss map}
All along this paper, and unless otherwise stated,
$\Sigma$ will denote a two-dimensional simply
connected smooth manifold endowed with a Riemannian flat metric
$\esiz,\esde$. We shall consider isometric
immersions $f:\Sigma\flecha\R^4$ with \emph{flat normal bundle}, that is
to say, isometric immersions whose normal curvature tensor vanishes
identically, $R^{\perp}\equiv 0$. It is then well known that there exists
a unit normal vector field $N\in\X^{\perp}(f)$ globally defined on
$\Sigma$ that is parallel in the normal bundle, i.e.
$dN(X)\in\X(f)$ for all $X\in \X(f)$ (see for instance \cite{Ten,Eno2}).
Any such unit normal vector field of $\X^{\perp}(f)$ will be called a
\emph{special section} on $\Sigma$.

It is immediate that special sections come in pairs. In fact, if
$N\in\X^{\perp}(f)$ is a special section, its orthonormal complement
$\hat{N}$ in $\X^{\perp}(f)$ is also a special section. Furthermore,
any other special section $N'$ on $\Sigma$ is given by
$N'=\cos\theta N+\sin\theta \hat{N}$ for some constant angle $\theta$.

On the other hand we shall denote by $N_1^f$ the \emph{first
normal space} of the immersion $f$ (see \cite{CaDa,DaTo}). As it was
explained in \cite{CaDa}, the study of isometric immersions from $\R^2$
into $\R^4$ needs a constancy condition on the dimension of $N_1^f$. Since
flat surfaces with ${\rm dim} (N_1^f)\equiv 1$ were completely described
in \cite{DaTo2}, we shall deal only with isometric immersions such that
${\rm dim} ( N_1^f )\equiv 2$. Rather than necessary, this condition is an
easy way to prevent {\it any} piece of the surface from lying in a
hyperplane of $\R^4$.

We plan to study flat surfaces with flat normal bundle in $\R^4$ by means of their
second fundamental form $\sigma$. For this, we will first of all show that if $f:\Sigma\flecha\R^4$ is
such a surface, we can endow $\Sigma$ with a natural Lorentz surface
structure, induced by $\sigma$. Let us recall (see \cite{Wns}) that a
\emph{Lorentz surface structure} on $\Sigma$ is the class of all
Lorentzian metrics on $\Sigma$ that are conformal to a specific Lorentzian
metric on $\Sigma$.

If $f:\Sigma\flecha\R^4$ is a flat surface with
$R^{\perp}\equiv 0$, and we choose local coordinates $(x,y)$ on
$\Sigma$ so that $\esiz,\esde=dx^2+dy^2$, and an orthonormal pair of
special sections $N,\hat{N}$ on $\Sigma$, then the following
\emph{structure equations} hold.
 \begin{equation}\label{structure}
 \left\{\def\arraystretch{1.2}
\begin{array}{l}
 f_{xx}=E_1N +E_2\hat{N} \\
 f_{yy}=G_1N +G_2\hat{N} \\
 f_{xy}= F_1N+F_2\hat{N}
 \end{array}\right. \hspace{0.2cm}\left\{\def\arraystretch{1.2} \begin{array}{ll}
 -N_x= E_1 f_x +F_1 f_y  \\
  -N_y= F_1 f_x +G_1 f_y
 \end{array}\right.\hspace{0.2cm} \left\{\def\arraystretch{1.2} \begin{array}{ll}
 -\hat{N}_x= E_2 f_x +F_2 f_y  \\
 -\hat{N}_y= F_2 f_x +G_2 f_y
 \end{array}\right.
 \end{equation}Here $E_i,F_i,G_i$ are smooth functions satisfying the
Gauss-Codazzi-Ricci equations, among which we quote the \emph{Gauss
equation} $E_1G_1-F_1^2 =-(E_2G_2-F_2^2)$, that is, $\det (dN)
+\det(d\hat{N})=0$.

Recall that the flat surface $f$ has ${\rm dim} (N_1^f)\equiv 2$. We claim
that for any $p\in\Sigma$ there exists a special section on $\Sigma$ which is an immersion over
a neighbourhood of $p$. Indeed, if this is not true, it follows from
\eqref{structure} and the fact that every special section $\xi$ is written
as $\xi=\cos\theta N+\sin\theta\hat{N}$ for some $\theta\in\R$ that
equations
 \begin{equation}\label{nevermore}
E_1G_1-F_1^2=0,\hspace{0.3cm} E_2G_2-F_2^2=0, \hspace{0,3cm} E_1G_2+E_2G_1
-2F_1F_2=0
  \end{equation}
must hold at $p$.

This implies the existence of $\landa_i,\mu_i\in\R$,
$i=1,2$, so that $\landa_i(E_i,F_i)+\mu_i(F_i,G_i)=0$ and $\landa_i,\mu_i$ do not
both vanish simultaneously. Besides, we have from ${\rm dim} (N_1^f)(p)=2$
and \eqref{nevermore} that $E_i$ and $G_i$ do not vanish at $p$. In this
way, by writing
$\nu_i=\mu_i/\landa_i\in\R$, it is obtained $$E_i=\nu_i^2 G_i, \hspace{0.3cm}
F_i=-\nu_iG_i, \hspace{0.3cm} i=1,2.$$ Finally, this last expression
together with the third equation in \eqref{nevermore} give
$\nu_1=\nu_2$, a contradiction with ${\rm dim} (N_1^f)(p)=2$. Summing
up, we have ensured the existence of a special section that is an immersion
in a neighbourhood of $p$ for every $p\in\Sigma$.

Now assume that $N$ is regular at $p$, and that $\{\xi,\hat{\xi}\}$ is another pair
of special sections on $\Sigma$. Then
there exists a smooth function $\landa$ on a neighbourhood of $p$
verifying $\esiz d\xi,d\hat{\xi}\esde=\landa \esiz dN,d\hat{N}\esde$,
that is, the Lorentzian pseudometric $\esiz d\xi,d\hat{\xi}\esde$
belongs to the Lorentz surface structure induced by $\esiz
dN,d\hat{N}\esde$. To see this, we first note the existence of vectors
$v_1,v_2\in T_p\Sigma$ that verify
 \begin{equation*}
 \left\{\def\arraystretch{1.2}
\begin{array}{l}
 \esiz dN_p(v_i),dN_p(v_i)\esde=1\\
 \esiz dN_p(v_1),dN_p(v_2)\esde=0 \\
 d\hat{N}_p(v_i)= k_i dN_p(v_i)
 \end{array}\right.
 \end{equation*} where $k_1,k_2\in\R$. Since $\xi=\cos\theta N+\sen \theta \hat{N}$ for
some $\theta$, we get that $\esiz d\xi_p(v_1),d\hat{\xi}_p(v_2)\esde=0=
\esiz dN_p(v_1),d\hat{N}_p(v_2)\esde$, and that $\esiz
d\xi,d\hat{\xi}\esde=\landa \esiz dN,d\hat{N}\esde$ holds at $p$ if and
only if the identity $$\sen\theta\cos\theta
(k_1^2k_2-k_2)=\sen\theta\cos\theta (k_1k_2^2 -k_1)$$ holds. But this is
true since the Gauss equation gives $\det (dN)+\det (d\hat{N})=0$, i.e.
$1+k_1k_2=0$.

\begin{remark}
Since at the regular points of $N$, its unit normal in the $3$-sphere is
$\hat{N}$, the relation $1+k_1k_2$ ensures that $N$ is a flat surface in
$\S^3$. In other words, every special section of $f$ defines, at its
regular points, a flat surface in $\S^3$. For any such surface, its second
fundamental form $\esiz dN,d\hat{N}\esde$ is a Lorentzian metric, and it
is classically known that the null directions of this metric point at the
asymptotic directions of the surface.
\end{remark}

To sum up, we have just proved the following result.
\begin{lem}
Let $f:\Sigma\flecha\R^4$ be a flat surface with flat normal bundle in
$\R^4$. There exists a natural Lorentz surface structure on $\Sigma$ such
that at every point $p$ its null directions are the asymptotic
directions of any of the special sections of $f$ that are regular at
$p$.
\end{lem}

Motivated by the above situation, next we introduce a generalization of
the flat surfaces in $\S^3$, which we call \emph{flat maps}, and for which
we allow the presence of singular points. These flat maps will be our
basic tool to construct flat surfaces in $\R^4$ with $R^{\perp}\equiv 0$.

\begin{defi}\label{fmap}
Let $\Sigma$ be a smooth simply connected surface.
A map $F:\Sigma\flecha\S^3$ on $\Sigma$ is called a \emph{flat map} if
there exist $\hat{F}:\Sigma\flecha\S^3$, $\omega:\Sigma\flecha\R$, and a canonical
immersion $$\Sigma\flecha u,v\text{-plane }$$ so that the following relations hold.
\begin{equation}\label{flatmap}
\def\arraystretch{1.5} \begin{array}{ll}
\esiz dF,dF\esde= du^2 +2\cos\omega \, dudv + dv^2, & \hspace{0.4cm}\esiz
F,\hat{F}\esde=0,
\\
\esiz dF,d\hat{F}\esde= 2\sin\omega \, dudv, & \hspace{0.4cm}\esiz
dF,\hat{F}\esde=\esiz F,d\hat{F}\esde=0,
\\
\esiz d\hat{F},d\hat{F}\esde= du^2 -2\cos\omega \, dudv + dv^2, &
\hspace{0.4cm} \omega_{uv}\equiv 0.
\\
 \end{array}
\end{equation}
\end{defi}

The immersion from $\Sigma$ into the $u,v$-plane must be thought of as a
coordinate immersion, but not as proper coordinates. This coordinate
immersion is inspired in the well known existence of a global Tschebysheff
net on every complete flat surface in $\S^3$ (see \cite{Spi}). We shall
refer to $\hat{F}$ (resp. $\omega$) as the \emph{polar map} (resp. the
\emph{angle}) of $F$. It is immediate that $\hat{F}$ is also a flat map
with polar map $-F$ and angle
$\hat{\omega}= \pi +\omega$. Moreover, if $F$ is an immersion, the
angle $\omega(u,v)$ can be chosen so that
$0<\omega<\pi$, and $F$ is a flat surface in $\S^3$ whose unit normal is $\hat{F}$.

At the end of the present Section we will show how to adapt the known
procedures for constructing complete flat surfaces in $\S^3$ to the
context of flat maps.

Next we turn our attention to the Gauss map, and its relation to the Lorentz surface structure
that we have just defined on $\Sigma$. Let $G_{2,4}$ denote the Grassmannian of
oriented $2$-planes in $\R^4$. The \emph{Gauss map} of the immersion $f:\Sigma\flecha\R^4$
is defined as the map $\cG:\Sigma\flecha G_{2,4}$ assigning to each
$p\in\Sigma$ the oriented tangent plane of $\Sigma$ at $p$.
Thus we have a Riemannian pseudometric on $\Sigma$ induced by $\cG$, and given by
$$\esiz ,\esde^*:=\esiz d\cG,d\cG\esde_{G_{2,4}},$$ where
$\esiz ,\esde_{G_{2,4}}$ is the Riemannian metric of $G_{2,4}$. In this way
$\cG:(\Sigma,\esiz,\esde^*)\flecha G_{2,4}$ becomes an isometric
immersion at its regular points.

Let us describe the pseudometric $\esiz,\esde^*$. For this we first
identify in the usual way $G_{2,4}$ with the complex quadric $Q_2$ of
$\C{\bf P}^3$ given in coordinates by $z_1^2+z_2^2+z_3^2+z_4^2=0$ (see \cite{HoOs}).
Then $Q_2$ is endowed with a Riemannian metric, the so-called
\emph{Fubini-Study} metric, given by
$$ds^2= 2\frac{ \sum_{j<k} |z_j dz_k -z_k dz_j|^2}{\left(\sum_{j=1}^4
|z_j|^2\right)^2}=2\frac{\cT (z,z)\cT
(dz,dz)-\left|\cT(z,dz)\right|^2}{\cT (z,z)^2},$$ where we are denoting by
$\cT$ the usual Hermitian product in $\C^4$, that is,
$$\cT(z,w)=\sum_{k=1}^4 z_k \overline{w}_k.$$ If we
choose local coordinates $(x,y)$ for the flat surface
$f:\Sigma\flecha\R^4$ such that $\esiz,\esde=dx^2+dy^2$, then $$\esiz ,\esde^*=
\cT(d(f_x+if_y),d(f_x+if_y))= \esiz d(f_x),d(f_x)\esde + \esiz
d(f_y),d(f_y)\esde.$$ From here and the structure equations
\eqref{structure} for $f$ we find that
 \begin{equation}\label{grass}
 \esiz,\esde^*=\esiz dN,dN\esde +\esiz d\hat{N},d\hat{N}\esde
 \end{equation}
where $N,\hat{N}$ is an arbitrary orthonormal pair of special sections on
$\Sigma$.

From now on, we shall assume that the Gauss map of the surface is regular, and thus that
$\esiz,\esde^*$ is a Riemannian metric on $\Sigma$. This is
not restrictive at all, since the following result holds, even in the case in which
$R^{\perp}\neq 0$ (see [Wei3]).
 \begin{lem}\label{1normal}
The Gauss map of $f$ is regular if and only if ${\rm dim} ( N_1^f )\equiv
2$ on $\Sigma$.
 \end{lem}
%

In this way, as it was shown above, the regularity of $\cG$ implies that for
every $p\in\Sigma$ there exists a
special section $N$ on $\Sigma$ that is an immersion over a neighbourhood
of $p$. Thus there is an open set $U\subset\Sigma$ containing $p$ so that
$N:U\flecha\S^3$ is a flat surface, and there exist local coordinates $s,t$ over
$U$ for which (see \cite{Spi} for instance)
$$\def\arraystretch{1.5} \begin{array}{ll}
\esiz dN,dN\esde&= ds^2 +2\cos\rho \, dsdt + dt^2
\\
\esiz dN,d\hat{N}\esde&= 2\sin\rho \, dsdt,
\\
\esiz d\hat{N},d\hat{N}\esde&= ds^2 -2\cos\rho \, dsdt + dt^2
\\
\end{array}$$ hold on $U$. Here $\rho:U\flecha\R$ satisfies $0<\rho<\pi$ and $\rho_{st}=0$,
and $\hat{N}$ is the orthogonal complement of the special section $N$. With this, it
follows directly from \eqref{grass} that
 \begin{equation}\label{cond1}
\esiz,\esde^*=2(ds^2+dt^2)
 \end{equation} on $U$. Particularly, $\cG:\Sigma\flecha G_{2,4}$ is a flat immersion.
In addition, since every orthonormal pair $(\xi,\hat{\xi})$
of special sections on $\Sigma$ is given by
$(\xi,\hat{\xi})=R_{\theta} (N,\hat{N})$, $R_{\theta}$ being a
rotation of angle $\theta$, then every such pair satisfies
 \begin{equation}\label{cond2} \esiz
\xi_s,\hat{\xi}_s\esde = \esiz \xi_t,\hat{\xi}_t\esde = 0
\end{equation} on $U$. In other words, $(\parc_s,\parc_t)$ point at the null
directions of the Lorentz surface $\Sigma$. It becomes clear then that
properties \eqref{cond1} and \eqref{cond2} characterize the parameters
$s,t$ on $U$. In this way, if $s',t'$ constitute a local
coordinate system satisfying both conditions on an open set
$U'$ in $\Sigma$ such that $U\cap U'\neq\emptyset$,
then $(\parc_s,\parc_t)= (\parc_{s'},\parc_{t'})$ over $U\cap U'$. From
here we are led to the main conclusion of this Section, that is stated in
terms of the notion of \emph{flat map}, which we introduced in Definition
\ref{fmap}.
 \begin{pro}\label{mainlema}
Let $f:\Sigma\flecha\R^4$ be a simply connected flat surface in
$\R^4$ with flat normal bundle and regular Gauss map $\cG$. There
exists an immersion from $\Sigma$ into the $u,v$-plane with respect
to which every special section $N$ on $\Sigma$ is a flat map whose
polar map
$\hat{N}$ is its orthonormal complement in $\X^{\perp}(f)$.
 \end{pro}
 \begin{proof}
Since $\cG$ is a flat immersion, $(\Sigma,\esiz,\esde^*)$ is a Riemannian
flat surface. Hence we can take a
$C^{\8}$ global conformal parameter $z:\Sigma\flecha\C$ such that $\esiz,\esde^*=
e^{2\varrho} |dz|^2$ for some smooth function $\varrho:\Sigma\flecha\R$.
The flatness of $\esiz,\esde^*$ implies that
$\varrho_{z\bar{z}}=0$, i.e. $\varrho$ is harmonic. Therefore there exists
a holomorphic map
$\phi:\Sigma\flecha\C$ such that $\varrho=\Re (\phi)$. Thus, if we take
the immersion from $\Sigma$ into the $u,v$-plane given by
$$\zeta=u+iv=\frac{\sqrt{2}}{2}\int e^{\phi} dz$$ it follows that
$$\esiz,\esde^*=|e^{\phi}dz|^2=2|d\zeta|^2 = 2\left(du^2 + dv^2\right).$$

Note that this last relation remains unchanged if we choose
$e^{i\varphi }\zeta$ instead of $\zeta$, where $\varphi\in\R$,
and that this change rotates
$(\parc_u,\parc_v)$ by an angle $\varphi$. Hence we can assume that at a
point $p\in\Sigma$, the above vectors point at the null directions of the
Lorentz surface $\Sigma$. But in that case, it follows from the existence and uniqueness
of the local parameters $(s,t)$ defined above that the coordinate immersion
$(u,v)$ satisfies properties \eqref{cond1} and \eqref{cond2} globally on
$\Sigma$.

In other words, we obtain an immersion of $\Sigma$ into $\R^2$ such that
\begin{align*}
 & \esiz,\esde^*= 2(du^2+dv^2) \\
 & \esiz N_u,\hat{N}_u\esde = \esiz N_v,\hat{N}_v\esde = 0
\end{align*}
hold on $\Sigma$, where here $(u,v)$ are canonical coordinates in $\R^2$
and $N,\hat{N}$ is an arbitrary orthonormal pair of special sections on
$\Sigma$. This is the canonical immersion appearing in the definition of flat map.

In this manner we have $\esiz N_u,N_u\esde +\esiz \hat{N}_u,\hat{N}_u\esde
=2$ and $$0=\left\esiz
\frac{(N+\hat{N})_u}{\sqrt{2}},\frac{(N-\hat{N})_u}{\sqrt{2}} \right\esde
= \frac{1}{2}\left(\esiz N_u,N_u\esde -\esiz
\hat{N}_u,\hat{N}_u\esde\right).$$ So, $\esiz N_u,N_u\esde =\esiz
\hat{N}_u,\hat{N}_u\esde =1$ and in the same way $\esiz N_v,N_v\esde
=\esiz \hat{N}_v,\hat{N}_v\esde =1$. Besides, from \eqref{cond1} we get
$\esiz N_u,N_v\esde=-\esiz\hat{N}_u,\hat{N}_v\esde$.

Summing up all these equations it is elementary to show the existence of a
smooth map
$\omega:\Sigma\flecha\R$ such that all conditions in \eqref{flatmap} are satisfied. We remark
that equation
$\omega_{uv}\equiv 0$ holds on the set $\cR$ of regular points of $N$,
since in there $N$ is a flat surface in $\S^3$, and it also holds in the
exterior of $\cR$, in where $\omega$ must be constant. Thus, by
continuity, we have $\omega_{uv}\equiv 0$ in all $\Sigma$.
 \end{proof}

Proposition \ref{mainlema} shows as a particular case that one can define
a global Tschebysheff coordinate immersion on every simply connected flat
surface in $\S^3$, regardless of its completeness. In other words, we have
proved the following
\begin{pro}
Every simply connected flat surface
in $\S^3$ (complete or not) is a flat map with $0<\omega<\pi$.
\end{pro}
To end up this Section, we describe two methods for constructing flat
maps, inspired in the known constructions of complete flat surfaces
in $\S^3$ given by Bianchi-Spivak and Kitagawa, respectively. These
classical constructions can be found in \cite{Spi,Kit1,Kit2,Wei2}
among others, and they carry over to the context of flat maps almost
unchanged, except for two details. One is that on a flat map the
parameters
$u,v$ are not proper coordinates, what produces some technical
difficulties. The other one is that on a flat map no regularity assumption
is required, and this simplifies the description of the construction
method. Hence, we will focus on these two differences, and just give a
sketch of the rest of the reasoning.

Let us begin by modifying the most classical construction, due to
Bianchi and put afterwards into a modern form by Spivak and Sasaki. Let us identify
$\R^4$ with the quaternions, so that $\S^2$ is the set of unit pure
quaternions, and $\S^3$ is that of unit quaternions.

Let $a_1(u),a_2(v)$ be two regular curves in $\S^3$ parametrized by
arclength, assume that
$a_1(0)=a_2(0)={\bf 1}$, and choose $\xi_0\in \S^3$ orthogonal to both $a_1'(0)$, $a_2'(0)$,
and lying in $T_{\bf 1}\S^3$. Then define $\xi_1(u)=a_1(u)\cdot
\xi_0$ and
$\xi_2(v)=\xi_0\cdot a_2(v)$, the dot denoting the product in $\S^3$.
Besides, consider the pair $(\Phi,\hat{\Phi})$
given by
\begin{equation}\label{Spivak}
\begin{array}{lll}
\Phi(u,v) &=& a_1(u)\cdot a_2(v), \\ \hat{\Phi}(u,v) &=& a_1(u)\cdot
\xi_0\cdot a_2(v),
\end{array}
\end{equation}
defined over a rectangle $R$ of the $u,v$-plane.

Finally consider an immersion $\Psi:\Sigma\flecha R$, where
$\Sigma$ is a simply connected surface, so that $\Psi(\Sigma)=R$,
and define the pair of maps $(F,\hat{F})$ from $\Sigma$ into $\S^3$ given
by
\begin{equation}\label{compo}
\left(F,\hat{F}\right)=\left(\Phi\circ\Psi, \hat{\Phi}\circ\Psi\right).
\end{equation}

Following in this general context \cite[Lemma 4.1]{Kit1} we obtain
\begin{lem}
If $\esiz a_i',\xi_i\esde\equiv 0$, $i=1,2$, then $F$ is a flat
map with polar map $\hat{F}$.
\end{lem}

Note that if we define the equivalence relation $\sim $ over $\Sigma$ so
that $p\sim q$ if and only if $\Psi(p)=\Psi(q)$, then it follows from
\eqref{compo} that $F,\hat{F}$ are well defined over $\Sigma^*=\Sigma
/\!\sim$. In addition, $\Psi:\Sigma^*\flecha R$ is well defined and
one-to-one, what shows that $u,v$ are proper coordinates on
$\Sigma^*$. Also note that $\Psi(\Sigma^*)=R$.

Next, we shall show that this construction can be reversed for any simply
connected analytic flat map, but not in general for smooth flat maps. To
do so, we begin by considering a particular situation, in which $F$ is a
flat map with polar map $\hat{F}$ and angle $\omega(u,v)$, defined on a simply connected
surface $\cU$, and endowed with the canonical immersion
$\Psi:\cU\flecha u,v\text{-plane}$, so that

1) $\Psi$ is one-to-one, and thus $u,v$ are proper coordinates on $\cU$,
and

2) $\Psi (\cU)$ is a rectangle
$R$ in the $u,v$-plane.

Then we can identify $\cU$ with $R$, and assume that
$(0,0)\in R$. We shall also suppose without losing generality that $F(0,0)={\bf 1}$ and
$\hat{F}(0,0)=\xi_0\in\S^3$. Now
define two regular curves in $\S^3$ as $a_1(u)=F(u,0)$ and
$a_2(v)=F(0,v)$. Let us compare, for the parameter $v$,
the curves in $\S^3$ given by $\Gamma_1(v)= F(u_0,v)$ and $\Gamma_2(v)=
a_1(u_0)\cdot a_2(v)$.

Consider the frame in $\S^3$ along $\Gamma_1$ given by
$$\{\mathbf{T}(v)=\Gamma_1'(v), \mathbf{N}(v)=\hat{F}_v (u_0,v),
\mathbf{B}(v)=\hat{F}(u_0,v))\}$$
and the corresponding frame along $\Gamma_2$ defined by $$\{\mathbf{t}(v)=a_1(u_0)\cdot a_2'(v), \mathbf{n}(v)=a_1(u_0)\cdot
\xi_0\cdot a_2'(v), \mathbf{b}(v)=a_1(u_0)\cdot \xi_0\cdot a_2(v)\}.$$

Since $(F,\hat{F})$ satisfy equations \eqref{flatmap}, we obtain that
$\nabla_{F_u} \hat{F} =F_u\times \hat{F}$, where $\nabla$ and $\times$ stand for the
Riemannian connection and the cross-product in $\S^3$, respectively. Thus
$\hat{F}$ is left invariant along $a_1(u)$, and it can be shown in the
same way that it is right invariant along $a_2(v)$. From there we obtain
that both frames coincide at $0$. Moreover, one can check that both frames
are the Frenet trihedron of their respective curves. Thus, by deriving we
find that the two curves have the same curvature and torsion in $\S^3$.
Since they have the same initial conditions, we conclude that
$\Gamma_1=\Gamma_2$. In other words, the flat map is recovered by means of
equations \eqref{Spivak} and \eqref{compo}.

Now assume that $F$ is an {\it analytic} flat map, defined on $\Sigma$.
Then its polar map $\hat{F}$ as well as the canonical immersion
$\Psi:\Sigma\flecha u,v$-plane are both analytic. Note that there exists
$\cU\subset\Sigma$ such that $\Psi|_{\cU}$ is one-to-one, and $\Psi(\cU)$
is a rectangle. Then, we have just proved above that \eqref{compo} holds
on $\cU$. Next, extend the analytic curves $a(u),b(v)$ as much as
possible, so that $\Phi,\hat{\Phi}$ are defined on a rectangle $R'$, and
they cannot be extended any further. Since all maps are analytic, we find
that \eqref{compo} holds not just in $\cU$, but in the largest
$W\subset\Sigma$ in which \eqref{compo} is well defined. But it
follows from the continuity of $F$ and the fact that $\Phi$ does not admit
a proper continuous extension away from $R'$ that $\Psi(\Sigma)\subset
R'$, and this indicates that $W=\Sigma$. Thus, \eqref{compo} holds on
$\Sigma$ globally. In particular $(F,\hat{F})$ are well defined over
$\Sigma^*=\Sigma / \sim$, and $u,v$ are proper coordinates. Note that
enlarging $\Sigma^*$ if necessary, the coordinates $u,v$ can be chosen to
be defined over all $R'$, so that $\Psi (\Sigma^*)=R'$. This completes the
first construction of flat maps.
\begin{remark}\label{coun}
There exist simply connected smooth flat maps that cannot be globally
recovered through equations \eqref{Spivak} and \eqref{compo}, even if
$\Psi$ is one-to-one. To show this, we begin by choosing three curves
$a_1,a_2,\tilde{a}_2:(-1,1)\flecha\S^3$ parametrized by arclength,
so that
\begin{enumerate}
\item
$a_1$ has torsion $\tau=1$, and $a_2,\tilde{a}_2$ have
$\tau=-1$.
\item
$a_1(0)=a_2(0)=\tilde{a}_2(0)={\bf 1}$
\item
$a_2(t)=\tilde{a}_2(t)$ for all $t\in (-1,0]$, but
$a_2(t)\neq \tilde{a}_2(t)$ for $t\in (0,1)$.
\end{enumerate}
Then define $\Sigma_0=(-1,1)\times(-1,1)\times \{0\}\subset \R^3$, and
$\Sigma_1=(-1,1)\times(-1,1)\times \{1\}\subset \R^3$, and consider
the map $\Phi:\Sigma_0\cup\Sigma_1\flecha \S^3$ given by

\begin{equation}\label{cha}
\begin{array}{lll}
\Phi(u,v,0)&=& a_1(u)\cdot a_2(v)\\ \Phi(u,v,1)&=& a_1(u)\cdot
\tilde{a}_2(v). \\
\end{array}
\end{equation}

It is clear that $\Phi$ is a non-connected flat map. Now define on
$\Sigma_0\cup\Sigma_1$ the equivalence relation that identifies
$(u,v,0) \simeq (u,v,1)$ for every $(u,v)\in (-1,1)\times (-1,0)$,
and consider the sets
$$\begin{array}{lll}
R_0 & = & \big((-1,1)\times (-1,0)\times \{0\}\big)\cup
\big((-1,-1/2)\times (-1,1)\times \{0\}\big)\subset \Sigma_0, \\ R_1 & = &
\big((-1,1)\times (-1,0)\times \{1\}\big)\cup \big( (1/2,1)\times
(-1,1)\times \{1\}\big) \subset \Sigma_1.
\end{array}$$ It then follows that $\Sigma:=(R_0\cup R_1)/\simeq$ is a simply
connected smooth surface, with proper coordinates $u,v$. In addition
$\Phi$ induces a well defined flat map $F$ on $\Sigma$ in the obvious way,
for which $u,v$ are Tschebysheff coordinates. However, $F$ cannot be
expressed by means of \eqref{Spivak} and \eqref{compo}, since it has two
non-congruent coordinate $v$-curves, namely, $a_2$ and $\tilde{a}_2$. Note
that this process works for smooth flat maps, but not for analytic flat
maps, since in this case it is impossible to reach the condition
$a_2(t)\neq \tilde{a}_2(t)$.
\end{remark}
%
Another difference between analytic flat maps and smooth flat maps is
presented in the next result.
\begin{cor}
Every simply connected analytic flat surface in $\S^3$ has globally
defined proper Tschebysheff coordinates $(u,v)$, and these can be extended
to be defined over a rectangle in the $u,v$-plane. In contrast, there
exist smooth simply connected flat surfaces in $\S^3$ that cannot be
endowed with global proper Tschebysheff coordinates.
\end{cor}
\begin{proof}
It only remains to show the assertion about the smooth case, i.e. we have
to construct a smooth flat map $F:\Sigma\flecha\S^3$ endowed with an
immersion $\Psi:\Sigma\flecha u,v$-plane such that $F$ is not well defined
over $\Sigma^*=\Sigma /\sim$, where $p\sim q$ if and only if $\Psi
(p)=\Psi(q)$.

For that purpose we shall add to the sets $R_0,R_1$ constructed in the
above Remark a new set on $\Sigma_0$, $$ \tilde{R}_0 =R_0\cup
\big((-1,1)\times (1/2,1)\times
 \{0\}\big)\subset\Sigma_0.$$
Now let $\Sigma:=\tilde{R}_0\cup R_1/\simeq$. Then $\Sigma$ is a simply
connected surface, and we can define the map
$\Psi:\Sigma\flecha u,v$-plane that takes each point $p$ of $\Sigma$ into its first
two coordinates, when $p$ is regarded as a point of
$\Sigma_0\cup\Sigma_1$. This map is an immersion, which is not one-to-one,
and whose image is not simply connected in the $u,v$-plane. Now consider
curves
$a_1,a_2,\tilde{a}_2:(-1,1)\flecha\S^3$ as in Remark \ref{coun},
and define
$\Phi(u,v):\Psi(\Sigma)\flecha\S^3$ as in \eqref{cha}.
This map is doubled valued exactly on the set of double points of the
immersion $\Psi$. Thus, $\Phi(u,v)$ may be regarded as a single valued map
$F:\Sigma\flecha\S^3$, which is obviously a flat map. Moreover,
$F$ is not well defined over $\Sigma^*$. This completes the
proof.
\end{proof}
As a consequence of this Corollary, it is not difficult to modify the
proof of \cite[Theorem 2.3]{Kit1} to obtain
\begin{cor}\label{noMo}
Every analytic flat surface in $\S^3$ (complete or not) is orientable.
\end{cor}

This result was proved for complete $C^{\8}$ flat surfaces in $\S^3$ by
Kitagawa in \cite{Kit1}. It contrasts with the situations of flat surfaces
in
$\R^3$ and in $\R{\bf P}^3$, where analytic flat Möbius strips are known
to exist (see \cite{ChKa,Wun}). As a consequence, we obtain that \emph{all
analytic flat Möbius strips in $\R{\bf P}^3$ come from analytic flat
cylinders in $\S^3$ that are invariant under the antipodal map}.

It is interesting to remark that the analyticity condition in Corollary
\ref{noMo} cannot be weakened to smoothness. Indeed, there are examples of
$C^{\8}$ embedded flat Möbius strips in $\S^3$, as we show next.

Let us define ${\rm Ad} (x) y= x\cdot y\cdot \bar{x}$, where $x,y\in\S^3$,
and recall the \emph{Hopf fibration} $h:\S^3\flecha\S^2$ given by $h(x)=
{\rm Ad}(x) \mathbf{i}$. Here the bar denotes conjugation in the
quaternions.

It was observed by H.B. Lawson that if $c$ is a regular curve in $\S^2$,
then
$h^{-1}(c)$ is a flat surface in $\S^3$. It has in general the topology of a
cylinder, but if $c$ is closed then $h^{-1}(c)$ is actually a torus.
Following \cite{Pin}, we shall call in general $h^{-1}(c)$ a {\it Hopf
cylinder}, and also a {\it Hopf torus} in case $c$ is closed.

\begin{figure}[h]\label{figur}
\mbox{}
\begin{center}
\includegraphics[clip,width=10cm]{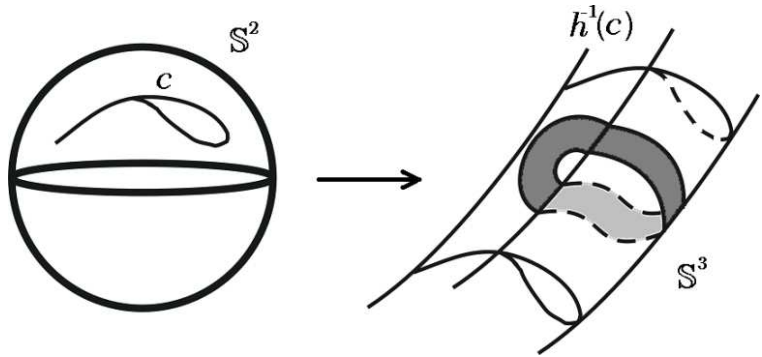}
\end{center}
\begin{quote}
\caption{A $C^{\8}$ flat Möbius strip embedded in $\S^3$, described as an
open set of a Hopf cylinder.}
\end{quote}
\end{figure}

With all of this, Figure $1$ shows an example of a $C^{\8}$ Möbius strip
lying on a Hopf cylinder in $\S^3$. For that, we choose a regular smooth
curve $c$ in $\S^2$ with the shape described in the figure (a curve with
this shape obviously cannot be real analytic). Then the shadowed region in
the Hopf cylinder $h^{-1}(c)$ clearly has the topology of a Möbius strip,
and it is trivially flat. This construction is inspired in the example of
a flat Möbius strip in $\R^3$ appearing in \cite{ChKa}.

Next, we turn to the second construction procedure of flat maps, which is
based on Kitagawa's work \cite{Kit1}. To explain it we shall use the
notation in \cite{Wei2}.

Let $U\S^2$ denote the bundle of unit tangent vectors to $\S^2$. Then we
can define the double cover
$\pi: \S^3\flecha U\S^2$ given by $$\pi(x)=\left({\rm Ad}(x) \mathbf{i},{\rm
Ad}(x) \xi_0\right),$$ where $\xi_0\in \S^3$ is orthogonal to
$\mathbf{1}$ and $\mathbf{i}$. Note that if $a(u)$ is a curve in $\S^3$,
and $c(u)$ denotes its Hopf projection in $\S^2$, it holds that
\begin{equation}\label{corch}
c'={\rm Ad} (a) [\bar{a}a', \mathbf{i}].
\end{equation}

Let $c_1(u),c_2(v)$ be two regular curves in $\S^2$ with $c_i(0)={\bf i}$,
$c_i'(0)=\xi_0$, and consider $\hat{c}_i=(c_i,c_i'/||c_i'||)$, with values
in $U\S^2$. Then there exist two regular curves $a_1(u),a_2(v)$ in $\S^3$
verifying $\pi (a_i)=\hat{c}_i$. That is, $c_i$ is the Hopf projection of
$a_i$, and so equation \eqref{corch} is verified, and besides $c_i'$ is
collinear with ${\rm Ad} (a_i) \xi_0$. These two conditions imply that
$$\esiz a_1', a_1\cdot \xi_0\esde\equiv 0 \equiv \esiz (\bar{a_2})',
\xi_0\cdot \bar{a_2}\esde.$$

At this point, the first procedure to construct flat maps ensures that the
map $F$ defined through equations \eqref{Spivak} (writing
$\bar{a_2}$ instead of $a_2$) and \eqref{compo} is a flat map, provided
that we parametrize first $a_i$ by arclength.

Conversely, suppose that we are given a flat map $F$ constructed by means
of \eqref{Spivak} and \eqref{compo}, and define
$\hat{c}_1(u),\hat{c}_2(v)$ two regular curves in $U\S^2$ as
$\pi(a_1)=\hat{c}_1$, $\pi(\bar{a_2})= \hat{c}_2$. Since $\esiz a_1',
a_1\cdot \xi_0\esde\equiv 0 \equiv \esiz a_2', \xi_0\cdot a_2\esde$, the
above computations ensure that if $c_1$ (resp. $c_2$) is the Hopf
projection of $a_1$ (resp. $\bar{a_2}$), then $c_1'$ (resp. $c_2'$) is
collinear with
${\rm Ad} (a_1) \xi_0$ (resp. ${\rm Ad} (\bar{a_2}) \xi_0$). In other
words, $\hat{c}_i=(c_i,c_i'/||c_i'||)$, what show that the above process
can be reversed. Summarizing, we have obtained the following
generalization of Kitagawa's construction.
\begin{teo}
Let $c_1(u),c_2(v)$ be two regular curves in $\S^2$, with
$c_i(0)={\bf i}$, $c_i'(0)=\xi_0$, for some $\xi_0\in\S^3$
orthogonal to both $\mathbf{1},\mathbf{i}$. Let $\pi :\S^3\flecha
U\S^2$ be the double cover given by
 \begin{equation}\label{doco}
\pi(x)=
\left({\rm Ad}(x) \mathbf{i},{\rm Ad}(x) \xi_0\right),
 \end{equation}
consider $a_1(u),a_2(v)$ two curves in $\S^3$ parametrized by arclength and
satisfying $\pi (a_i)=(c_i,c_i'/||c_i'||)$, and define
\begin{equation}\label{Kita}
\begin{array}{lll}
\Phi(u,v) &=&  a_1(u)\cdot \bar{a_2}(v) \\
\hat{\Phi}(u,v) &=& a_1(u)\cdot \xi_0\cdot \bar{a_2}(v).
\end{array}
\end{equation}
on a rectangle $R$ in the $u,v$-plane. If $\Sigma$ is a simply connected
surface and
$\Psi:\Sigma\flecha \Psi(\Sigma)=R$ is an immersion, then the pair
$(F,\hat{F})$ given by
\begin{equation*}
\left(F,\hat{F}\right)=\left(\Phi\circ\Psi, \hat{\Phi}\circ\Psi\right).
\end{equation*}
constitutes a flat map $F$ with polar map $\hat{F}$. Conversely, every
analytic flat map is constructed in this way for some $\xi_0$.
\end{teo}
\begin{remark}
In the known constructions of flat surfaces in $\S^3$ that we have just
generalized, one has to assume further conditions to forbid the appearance
of singular points. It is not difficult to check that if one assumes those
additional conditions in the above descriptions of flat maps, it is
obtained a flat surface in $\S^3$. It is remarkable that this general
global process works for \emph{any} analytic simply connected flat surface
in $\S^3$, and not just locally or under a completeness assumption.

In addition, arguing as in Remark \ref{coun} we can produce $C^{\8}$
simply connected flat surfaces in $\S^3$ that cannot be globally obtained
through any of these two methods. Indeed, these examples will posses three
mutually non-congruent asymptotic curves.
\end{remark}

\section{Construction of flat surfaces}
Let $f:\Sigma\flecha\R^4$ be an isometric immersion of a simply connected
flat surface $\Sigma$ into $\R^4$ with flat normal bundle and regular
Gauss map $\cG$. Let $N:\Sigma\flecha\S^3$ be a special section on
$\Sigma$. Then, from Lemma \ref{mainlema} we find that $N$ is a flat map. Besides,
if we consider the immersion of $\Sigma$ into the $u,v$-plane constructed
in that Lemma, it is obtained that
$\{N,\hat{N},N_u,\hat{N}_u\}$ is an orthonormal frame of $\R^4$. Here
$\hat{N}$ stands for the polar map of $N$, which is also a special section
on $\Sigma$. From here the isometric immersion
$f$ can be expressed as $$f=\alfa N +\beta\hat{N} +\alfa_u N_u +\beta_u \hat{N}_u,$$
where $\alfa=\esiz f,N\esde$ and $\beta=\esiz f,\hat{N}\esde$. Noting that
 \begin{equation}\label{otrosis}
\left(\def\arraystretch{1.0}\begin{array}{c}
 N_v \\ \hat{N}_v\\
\end{array}\right)= \left(\def\arraystretch{1.1}\begin{array}{cc}
 \cos\omega & \sin \omega \\ \sin\omega &-\cos\omega\\
\end{array}\right)\left(\def\arraystretch{1.0}\begin{array}{c}
 N_u \\ \hat{N}_u\\
\end{array}\right),
 \end{equation}
we obtain that $\alfa,\beta$ satisfy the differential system
 \begin{equation*}
\left(\def\arraystretch{1.0}\begin{array}{c}
 \alfa_v \\ \beta_v\\
\end{array}\right)= \left(\def\arraystretch{1.1}\begin{array}{cc}
 \cos\omega & \sin \omega \\ \sin\omega &-\cos\omega\\
\end{array}\right)\left(\def\arraystretch{1.0}\begin{array}{c}
 \alfa_u \\ \beta_u\\
\end{array}\right).
 \end{equation*}

This process can be reversed. More concretely, we get the following
result, which describes all flat surfaces in $\R^4$ with flat normal
bundle and regular Gauss map in terms of flat maps.
 \begin{teo}\label{repre}
Let $f:\Sigma\flecha\R^4$ be an isometric immersion of a simply connected
flat surface $\Sigma$ into $\R^4$ with flat normal bundle and regular
Gauss map. There exists an immersion of $\Sigma$ into the $u,v$-plane with
respect to which every special section $N$ on $\Sigma$ is a flat map, and
if $\hat{N}$ is the polar map of any such $N$ the immersion $f$ is
expressed as
 \begin{equation}\label{rep}
 f=\alfa N +\beta\hat{N} +\alfa_u N_u +\beta_u \hat{N}_u.
 \end{equation}
Besides, if $\omega$ denotes the angle of the flat map $N$, the functions
$\alfa,\beta:\Sigma\flecha\R$ are solution of the differential system
 \begin{equation}\label{sistema}
\left(\def\arraystretch{1.0}\begin{array}{c}
 \alfa_v \\ \beta_v\\
\end{array}\right)= \left(\def\arraystretch{1.1}\begin{array}{cc}
 \cos\omega & \sin \omega \\ \sin\omega &-\cos\omega\\
\end{array}\right)\left(\def\arraystretch{1.0}\begin{array}{c}
 \alfa_u \\ \beta_u\\
\end{array}\right).
 \end{equation}
Conversely, let $N:\Sigma\flecha\S^3$ be a flat map with angle $\omega$
and polar map $\hat{N}$, and let $\alfa,\beta:\Sigma\flecha\R$ be a
solution of the differential system \eqref{sistema}. Then the mapping
$f:\Sigma\flecha\R^4$ given by \eqref{rep} constitutes at its regular points
a flat surface in $\R^4$ with flat normal bundle and regular Gauss map,
for which $N,\hat{N}$ are orthonormal special sections.
 \end{teo}
\begin{remark}
By considering its universal covering if necessary, every flat surface in
$\R^4$ with $R^{\perp}\equiv 0$ and regular Gauss map is globally described via
Theorem {\rm \ref{repre}}. Conversely, a flat surface
$f:\Sigma\flecha\R^4$ produced by \eqref{rep} may actually be well defined
on some quotient of
$\Sigma$ with non-trivial topology. We shall follow this last strategy in our
construction of flat tori and flat cylinders in the next Section.
\end{remark}

Let us examine more closely this representation formula.

To begin with, it is convenient to characterize the singular points of the
flat surface
$f$. Actually, the geometric meaning of the system \eqref{sistema}
will provide a geometric interpretation of the appearance of singular
points in $f$, as we show next.

First of all, note that \eqref{sistema} gives
$$\def\arraystretch{1.3}\begin{array}{lll} \alfa_{u v}&=& \cos\omega
(\alfa_{uu} +\omega_u\beta_u) +\sin \omega (\beta_{uu} -\omega_u \alfa_u),
\\ \beta_{u v}&=& \cos\omega (-\beta_{uu} +\omega_u\alfa_u) +\sin \omega
(\alfa_{uu} +\omega_u \beta_u). \\
\end{array}$$ If we denote

\begin{equation}\label{ambm}
A=\alfa +\alfa_{uu} +\omega_u\beta_u,\hspace{0.2cm} B=\beta +\beta_{uu}
-\omega_u\alfa_u,
\end{equation} and
 \begin{equation*}
\left(\def\arraystretch{1.0}\begin{array}{c}
 \hat{A} \\ \hat{B}\\
\end{array}\right)= \left(\def\arraystretch{1.1}\begin{array}{cc}
 \cos\omega & \sin \omega \\ \sin\omega &-\cos\omega\\
\end{array}\right)\left(\def\arraystretch{1.0}\begin{array}{c}
 A \\ B\\
\end{array}\right),
 \end{equation*}
then we obtain
 \begin{equation}\label{otrama}
\left(\def\arraystretch{1.0}\begin{array}{c}
 f_u \\ f_v\\
\end{array}\right)= \left(\def\arraystretch{1.1}\begin{array}{cc}
 A & B \\ \hat{A} & \hat{B}\\
\end{array}\right)\left(\def\arraystretch{1.0}\begin{array}{c}
 N_u \\ \hat{N}_u\\
\end{array}\right).
 \end{equation}
That is, if $x_1,x_2$ are canonical coordinates in the
$\{N_u,\hat{N}_u\}$-plane at some point $p\in \Sigma$, then the vectors
$f_u,f_v$ are symmetric with respect to the axis $\sin(\omega
/2)x_1=\cos(\omega /2)x_2$. Thus the surface
$f:\Sigma\flecha\R^4$ has a singular point at $p$ if and only if $f_u(p)$
lies on the above axis or is perpendicular to it. That is, $f$ is regular
if and only if the relation
 \begin{equation}\label{sing}
 (A^2-B^2)\sin \omega -2AB\cos \omega \neq 0
 \end{equation}
holds at every point. Moreover, the metric of $f$ in the coordinates
$(u,v)$ is written as
 \begin{equation}\label{metrica}
 \esiz df,df\esde=(A^2+B^2)(du^2+dv^2) +2 \left((A^2-B^2)\cos\omega +2AB
\sin\omega\right) dudv.
 \end{equation}
\begin{remark}\label{abab}
If $\alfa,\beta$ are solutions of \eqref{sistema}, then the functions
$A,B$ appearing in \eqref{ambm} are also solutions of \eqref{sistema}.
This follows from the condition $(f_u)_v=(f_v)_u$ in \eqref{otrama},
taking into account that the second derivatives of $N$ and $\hat{N}$ can
be expressed in terms of $\{N,\hat{N},N_u,\hat{N}_u\} $ by differentiation
of the system \eqref{otrosis}.
\end{remark}

In the remaining of this Section we shall discuss how to construct
solutions of the differential system \eqref{sistema}.

In the particular case in which the angle $\omega$ is constant, the system
can be completely integrated, and $(\alfa,\beta)$ are of the form
$(\alfa,\beta)(u,v)=(\alfa_1,\beta_1)(u+v)+(\alfa_2,\beta_2)(u-v)$ for arbitrary real
functions $\alfa_i,\beta_i$. Geometrically,
the flat surfaces in $\R^4$ for which $\omega$ is constant are
precisely the products of curves in $\R^4$, that is, the products
$\gamma_1\times \gamma_2$ of two regular plane curves in $\R^4$ lying in orthogonal
planes. Those are the flat surfaces whose Gauss image in $G_{2,4}\equiv
\S^2\times \S^2$ lies in the product of two great circles.

On the other hand, the system \eqref{sistema} admits constant solutions.
If $\alfa,\beta $ are constant, then $f(\Sigma)$ lies in a
$3$-sphere in $\R^4$ centred at the origin. The situation in which the
image of $f$ lies in an affine $3$-sphere not centred at ${\bf 0}$ is more
interesting from the viewpoint of the system \eqref{sistema}, and will be
crucial in the construction of new flat tori in $\R^4$ that we will
accomplish in the next Section. If $f$ is a flat surface in a
$3$-sphere of radius $\rho$ and centred at $a\in\R^4$, then there is some
special normal section $N$ of $f$ such that $f= a +\rho N$. If
$\hat{N}$ denotes the polar map of the flat map $N$, it follows that $\alfa = \esiz a,N\esde +
\rho$, and $\beta= \esiz a,\hat{N}\esde$. Thus, we conclude from the
linearity of \eqref{sistema} the following.

\begin{pro}\label{geo}
Let $N$ be a flat map with angle $\omega$, and denote its polar map
by $\hat{N}$. Then the coordinates $(N_i,\hat{N}_i)$, $1\leq i\leq
4$, are solutions of the differential system \eqref{sistema}.
\end{pro}

These solutions, as well as their linear combinations, will be referred to
as \emph{geometric solutions} of the system \eqref{sistema}. The above
comments ensure that {\it a flat surface in $\R^4$ given through equation
\eqref{rep} lies in a $3$-sphere of $\R^4$ if and only if $\alfa,\beta$
are geometric solutions of \eqref{sistema} (or constants)}.

It is important to observe that in general this geometric solutions cannot
be obtained explicitly in analytic terms. For instance, the coordinates of
any Hopf torus $h^{-1}(c)$ in $\S^3$ have an explicit formula in terms of
the curve $c$ in $\S^2$ that generates it. However, it is not generally
possible to derive explicitly the expression of these coordinates with
respect to the parameters
$(u,v)$.

Apart from these geometric integrations, the system \eqref{sistema} also
admits some particular analytic integrations, as we show next.

First, let us fix some matrix notation to treat \eqref{sistema}. For that,
we define ${\bf X}=(\alfa,\beta)^T$, and ${\bf M}$ the matrix such that
\eqref{sistema} is ${\bf X}_v={\bf M}{\bf X}_u$. Also note that, since
$\omega_{uv}\equiv 0$, the angle $\omega$ can be written as
$\omega(u,v)=\omega_1(u)+\omega_2(v)$. Thus ${\bf M}(u,v)={\bf L}(u) {\bf
H}(v)$, where
\begin{equation*}
{\bf L}(u)= \left(\def\arraystretch{1.1}\begin{array}{cc}
 \cos\omega_1 & \sin \omega_1 \\ \sin\omega_1 &-\cos\omega_1\\
\end{array}\right), \hspace{0.3cm}
{\bf H}(v)= \left(\def\arraystretch{1.1}\begin{array}{cc}
 \cos\omega_2 & \sin \omega_2 \\ -\sin\omega_2 &\cos\omega_2\\
\end{array}\right).
 \end{equation*}

We remark that if ${\bf X}$ is a solution of \eqref{sistema}, then ${\bf Z}={\bf
M X}_{uv}$ is another solution of \eqref{sistema}.

This fact can be reversed to obtain a more interesting situation. Let us
start with a solution ${\bf X}$ of \eqref{sistema}, so that ${\bf X}$
verifies ${\bf X}_v= {\bf L}(u) {\bf H}(v){\bf X}_u$. Then $({\bf L
X})_v=({\bf H X})_u$, what ensures the existence of a vector ${\bf Y}$
such that ${\bf Y}_u ={\bf L X}$ and ${\bf Y}_v ={\bf H X}$. Thus $({\bf
H^{-1} Y})_u=({\bf L Y})_v$, and there exists some ${\bf Z}$ verifying
${\bf Z}_u ={\bf L Y}$ and ${\bf Z}_v ={\bf H^{-1} Y}$. But this finally
implies that ${\bf Z}_v= {\bf M Z}_u$, that is, ${\bf Z}$ is again a
solution of \eqref{sistema}. To sum up, the vector $${\bf Z}(u,v)=\int
{\bf LY}du +\int {\bf H^{-1} Y} dv,$$ where ${\bf Y}=\int {\bf L X} du
+\int {\bf H X} dv$ and ${\bf X}$ is an arbitrary solution of
\eqref{sistema}, is again a solution of \eqref{sistema}.

For instance, if we let ${\bf X}=(0,0)$, which is a trivial solution of
\eqref{sistema}, and ${\bf v_0}$ is a fixed vector of $\R^2$,
then $${\bf Z}(u,v)=\left(\int \left(\def\arraystretch{1.1}\begin{array}{cc}
 \cos\omega_1 & \sin \omega_1 \\ \sin\omega_1 &-\cos\omega_1\\
\end{array}\right) du  + \int \left(\def\arraystretch{1.1}\begin{array}{cc}
 \cos\omega_2 & \sin \omega_2 \\ -\sin\omega_2 &\cos\omega_2\\
\end{array}\right)dv\right) {\bf v_0}$$ is a new non-constant
solution of \eqref{sistema}. In particular, choices like
$$\left\{\def\arraystretch{1.3}\begin{array}{lll}
\alfa_1&=&\int \cos \omega_1 du +\int\cos\omega_2 dv,\\
\beta_1&=&\int \sin \omega_1 du +\int\sin\omega_2 dv,\\
\end{array}\right. \hspace{0.3cm}  \left\{\def\arraystretch{1.3}\begin{array}{lll}
\alfa_2&=&\int \sin \omega_1 du -\int\sin\omega_2 dv,\\
\beta_2&=&-\int \cos \omega_1 du +\int\cos\omega_2 dv,\\
\end{array}\right.$$ constitute solutions of \eqref{sistema}. These
solutions have the particular feature that $(\alfa_i,\beta_i)=(A_i,B_i)$,
where $A_i,B_i$ are defined as in Remark \ref{abab}, for $i=1,2$.

Next, we shall consider a special case in which all solutions of
\eqref{sistema} can be found. As it was shown in the preceding Section,
every flat map can be constructed as the product of two curves in $\S^3$
with some special properties, in analogy with the classical Bianchi-Spivak
method. Since the choice of circles in this process only gives Clifford
tori, the most simple situation that is non-trivial is to consider helices
in
$\S^3$ with torsion $1$ and $-1$, respectively.

Let us consider the curve $\sigma(s):\R\flecha\S^3$ given by
 \begin{equation}\label{heli}
 \sigma(s)= \frac{1}{\sqrt{1+r^2}}\left(r\cos(s/r),r\sin (s/r),\cos (r s),
 \sin (r s)\right),\hspace{0.2cm} r>1.
 \end{equation}

Straightforward computations show that $\sigma(s)$ is parametrized by
arclength, its curvature is constant, $\kappa\equiv (r^2-1)/r$, and its
torsion $\tau$ verifies $\tau^2=1$. Thus, it is a helix in $\S^3$.
Conversely, every helix in $\S^3$ with $\tau^2=1$ is written in that way,
modulo a (possibly orientation-reversing) isometry of $\S^3$. Moreover,
$\sigma(s)$ is periodic if and only if $r^2\in \Q$. In the case where
$r\in \N$, $\sigma(s)$ is a simple closed curve of period $2\pi r$.

Now consider two of these helices, one with $\tau=1$, the other with
$\tau=-1$, and both with the same curvature $\kappa=(r^2-1)/r$, and
construct a flat map $N:\R^2\flecha\S^3$ following the Bianchi-Spivak
method. If we let $\mu=(r^2-1)/2r$, the angle of this flat map is
$\omega(u,v)= 2\mu (u+v)$. Let $\alfa,\beta$ be solutions of
\eqref{sistema} for this angle, and define $$\phi(u,v)=\cos
\left(\frac{\omega}{2}\right)\alfa +\sin \left(\frac{\omega}{2}\right)
\beta .$$ Then $\phi_u=\phi_v$, and thus $\phi(u,v)=2g'(u+v)$ for some
smooth real function $g(t)$. From here, if we define $\psi(u,v)= 2\mu
g(u+v) +h(u-v)$, $h(t)$ being a smooth real function, the system can be
completely integrated, and $\alfa,\beta$ are given by

 \begin{equation}\label{inhe}
\left(\def\arraystretch{1.0}\begin{array}{c}
 \alfa \\ \beta\\
\end{array}\right)= \left(\def\arraystretch{1.1}\begin{array}{cc}
 \phi(u,v) & \psi(u,v) \\ -\psi(u,v) & \phi(u,v))\\
\end{array}\right)\left(\def\arraystretch{1.0}\begin{array}{c}
\cos\left(\mu(u+v)\right)\\\sin\left(\mu(u+v)\right)\\
\end{array}\right).
\end{equation}
In this way we obtain via \eqref{rep} the explicit coordinates of
infinitely many flat surfaces $f:\Sigma\flecha\R^4$, possibly with
singular points, with flat normal bundle and regular Gauss map. An
analogous integration is obtained if we assume that the angle is
$\omega(u,v)=\mu (u-v)$.

The regularity condition for these surfaces is written as
$$\left(2G'(u+v)\right)^2 +\left(2\mu G(u+v) +H(u-v) \right)^2 \neq 0,$$
where $G(t)= (1+\mu^2) g(t) +g''(t)$ and $H(t)= (1+\mu^2) h(t) + h''(t)$.
With this, it is easy to arrange $g,h$ in order to define $f$ without
singular points. Actually, one can construct $f$ to be defined on the
whole $u,v$-plane. Furthermore, if we choose $g(t)=0$ and $\mu$ so that
$\mu^2\in \Q$, then $f:u,v\text{-plane}\flecha\R^4$ is a cylinder.

Any surface of this kind is characterized by the fact that its Gauss image
in $G_{2,4}\equiv \S^2\times\S^2$ lies in the product of two circles with
geodesic curvatures verifying $k_1^2=k_2^2$. In the case where the domain
of the parameters $(u,v)$ is the whole $u,v$-plane, its Gauss map is
compact, even though the surface $f$ is not if $\omega\neq \text{ const}$.
This contrasts with the situation in $\S^3$, where all flat surfaces with
globally defined Tschebysheff coordinates and whose Gauss map is a compact
surface in $\S^2\times\S^2$ are actually flat tori \cite{Wei1}.

This special shape of the Gauss map implies the following facts:
\begin{enumerate}
\item
None of these flat surfaces $f:\Sigma\flecha\R^4$ has a special section
without singular points. This situation justifies the way in which we
introduced in Section $2$ the Lorentz surface structure on $\Sigma$, since
$\esiz dN,d\hat{N}\esde$ is not in general a non-degenerate Lorentzian
metric for any pair $(N,\hat{N})$ of special sections.
\item
These flat surfaces are exactly the ones made up by compositions (see
\cite{DaTo}). This follows from the characterization of compositions in
terms of their Gauss maps that was observed in \cite{Wei3}.
 \item
In particular, as it was also observed in \cite{Wei3}, none of these flat
surfaces is complete (unless $\mu=0$).
\end{enumerate}

To end up this Section, we shall construct explicitly some new flat
cylinders in $\R^4$ with $R^{\perp}\equiv 0$ and regular Gauss map that
are not compositions.

To do so, we multiply again two helices in $\S^3$, this time with
different curvatures, so that the angle of the flat map that they
determine is $\omega(u,v)=2ru+2sv$, where $r,s\in\R$, $r^2\neq s^2$. It is
straightforward to check that $$
\left\{\def\arraystretch{1.2}\begin{array}{lll} \alfa & = & \exp (su+rv)
\big( \cos (ru+sv) +\sin (ru+sv) \big)\\ \beta & = & \exp (su+rv) \big(
\sin (ru+sv) -\cos (ru+sv)\big)
 \end{array}\right.$$
are solutions of \eqref{sistema} for this angle. Since they satisfy the
regularity condition \eqref{sing}, they give rise to a flat surface with
$R^{\perp}\equiv 0$ and regular Gauss map. If $r,s\in\N$, these two
helices are closed and the resulting immersion $f:u,v\text{-plane}\flecha
\R^4$ is a cylinder. However, none of this examples is complete. We shall
show at the end of the paper how to produce new complete flat cylinders in
$\R^4$.

\section{New flat tori}
In this Section we shall exhibit a general family of new flat tori in
$\R^4$. This family, which is the first one that comes out since Bianchi's
works in the 19th century on flat tori in
$\S^3$, is made up by flat tori with flat normal bundle, and provides a
counterexample to Borisenko's conjecture (see
\cite{Bor}). Specifically, we shall give a constructive proof of the following
result.
\begin{teo}\label{Boris}
There exists a general family of flat tori in $\R^4$ with flat normal
bundle that are not product of curves, and do not lie in any affine
$3$-sphere of
$\R^4$.
\end{teo}

It is remarkable that the tori of the family that we construct cannot be
given in explicit coordinates.

The key step of the proof consists in verifying the following fact, that
has interest in its own right.
\begin{teo}\label{truco}
Given $n\in \N$, $n>1$, there exists a Hopf torus $T$ in $ \S^3 $
with non-constant angle
$\omega(u)$ such that the Hopf cylinder in $\S^3$ with angle
$\widetilde{\omega}(u)=\omega(u/n)$ is actually a torus. Moreover,
$T$ can be chosen arbitrarily close to any prescribed Clifford torus in
$\S^3$.
\end{teo}

We shall actually show that there are infinitely many Hopf tori in
$\S^3$ with the above properties for any $n\in \N$, $n>1$. Here, by a
Clifford torus we mean a product torus $\S^1(r)\times \S^1(\sqrt{1-r^2})$
in
$\S^3$, which can be constructed as the lift via the Hopf fibration $h$ of
a circle in $\S^2$. The \emph{closeness} stated in Theorem \ref{truco}
must be understood in the following sense: the curve $c\in \S^2$ that
generates the Hopf torus
$T$ can be chosen arbitrarily close to a circle in $\S^2$.
\begin{remark}
Let $c(u)$ be a closed regular curve in $\S^2$, the parameter $u$ being
the arclength parameter of its asymptotic lift in $\S^3$ (see {\rm
\cite{Kit1}}). If $k(u)$ denotes the geodesic curvature of $c(u)$, the
angle $\omega(u)$ of the Hopf torus $h^{-1}(c)$ is given by
$$\omega(u)=\cot^{-1} (k(u))$$ Thus Theorem \ref{truco} asserts the
existence of a periodic curve $c(u)$ in
$\S^2$ such that the only (up to congruence) curve $\tilde{c}(u)$ in
$\S^2$ with geodesic curvature
$k(u/n)$ is periodic. This is strongly related to the
following problem posed by S.S. Chern (see {\rm \cite{CC}}): when is a
curve with periodic curvatures in a space form periodic?
\end{remark}
\begin{remark}\label{alp}
It was proved in {\rm \cite{Kit1}} that the arclength parameter $s$ of a
curve
$c(s)$ in $\S^2$ is expressed in terms of the arclength parameter $u$ of
its asymptotic lift in $\S^3$ as
 \begin{equation}\label{relpar}
 \frac{ds}{du}= \frac{2}{\sqrt{1+k(u)^2}},
 \end{equation}
where $k(u)$ is the geodesic curvature of the curve $c(s)$ at $u=u(s)$.
\end{remark}
\begin{proof1}
Fix $n\in \N$, $n>1$. From Theorem \ref{truco} we know the existence of a
Hopf torus $N$ in $\S^3$ with non-constant angle $\omega (u)$ such that
$\tilde{\omega}(u)=\omega(u/n)$ is also the
angle of a Hopf torus $\tilde{N}$ in $\S^3$. Let
$\tilde{\alfa},\tilde{\beta}$ be geometric solutions corresponding to
$\tilde{\omega}$. Then they are of the form
$(\tilde{\alfa},\tilde{\beta})= {\bf \tilde{a}}(u)\sin v +{\bf \tilde{b}} (u)\cos
v$. Now let us define
 \begin{equation}\label{ab}
(\alfa,\beta)(u,v)=(\tilde{\alfa},\tilde{\beta}) (nu,nv).
 \end{equation}
It is easy
to check that $\alfa,\beta $ are solutions of \eqref{sistema} for the
angle $\omega$ we started with. Moreover, they are not geometric
solutions, since they are not of the form
$(\overline{\alfa},\overline{\beta})= {\bf a}(u)\sin v +{\bf b} (u)\cos
v$.

Now, since $N$ is a Hopf torus in $\S^3$, the map $N(u,v)$ is doubly periodic in $u$ and
$v$. More exactly, there exists $T >0$ such that $N(u,v)$ is well
defined on the rectangular torus $\R^2 / \Landa$, where $\Landa=\{(T p,2\pi q)\in\R^2 :
p,q\in \Z\}$. Thus $\omega(u)$ is $T$-periodic,
and $\tilde{\omega}(u)$ is $nT$-periodic. In this way
the Hopf torus $\tilde{N}$ is defined on $\R^2 / \Gamma$, where
$\Gamma=\{(T k n p ,2\pi q) : p,q\in\Z\}$, for some $k\in \N$. Thus $\alfa,\beta$ and
$N$ are defined on $\R^2 /\Gamma $.

Summing up, the map $f$ given by \eqref{rep} is well defined
on a rectangular torus $\T^2=\R^2 /\Gamma$ with respect to the coordinates
$(u,v)$. Note that it is not contained in any affine $3$-sphere,
since $\alfa,\beta$ are not geometric solutions of \eqref{sistema}.
In other words, we have obtained a flat torus in $\R^4$,
possibly with singular points, that has flat normal bundle and regular Gauss map.
Furthermore, it is not a product of curves and it does not lie in any $3$-sphere of $\R^4$.

To conclude the proof of the Theorem we need to show that this flat torus
can be chosen without singular points. For this purpose, we begin by
defining on $\T^2$, for any $\landa >0$, the functions $\alfa_{\landa} = 1+\landa
\alfa$ and $\beta_{\landa}=\landa\beta$. It follows that
$\alfa_{\landa},\beta_{\landa}$ also constitute a solution of \eqref{sistema},
and thus they give rise to a new flat surface in $\R^4$, $f_{\landa}:\T^2\flecha\R^4$,
possibly with singular points. Let us also note that
$f$ lies in an affine $3$-sphere if and only if so does $f_{\landa}$. Now
observe that $(\alfa_{\landa},\beta_{\landa})$ tends pointwise to $(\alfa_0,\beta_0):=(1,0)$
as $\landa$ tends to $0$. But since $$|\alfa_{\landa} (u,v)-\alfa_{\landa'}
(u,v)|=|\landa -\landa'| |\alfa (u,v)|,$$ and the same relation holds for
$\beta_{\landa}$, we infer that $(\alfa_{\landa},\beta_{\landa})$ tends
uniformly on $\T^2$ to $(\alfa_0,\beta_0)$. Since $(\alfa_0,\beta_0)$ were
chosen to satisfy the regularity condition \eqref{sing}, we obtain the
existence of some
$\landa_0 >0$ such that $f_{\landa}$ is an immersion of $\T^2$ into $\R^4$
for any $\landa <\landa_0$. This completes the proof.
\end{proof1}
\begin{proof2}
Let us begin by defining $\cC$ as the set of smooth closed regular curves
in $\S^2$, parametrized as $\alfa (t):[0,1]\flecha\S^2$, where $t=s/\ell$,
$s$ denotes the arclength parameter, and $\ell$ is the length of the
curve. On $\cC$ we can define a topology through the norm
$||\cdot ||_{\cC}$ given by $$|| \alfa (t) ||_{\cC}=
||\alfa (t)||_{\8} + ||\alfa' (t)||_{\8}+||\alfa'' (t)||_{\8}.$$ The key
ingredient of this proof is the construction for each $n\in\N$, $n>1$, of
a continuous map
$\cA_n$ from $\cC$ into $\R$ with
the following property {\bf (P)}: if $\omega(u)$ denotes the angle of the Hopf torus
$h^{-1}(\alfa(t))$, then the Hopf cylinder in $\S^3$ of angle
$\omega(u/n)$ is actually a torus if and only if $\cA_n (\alfa(t))\in \Q$.

Fix $n\in\N$, $n>1$. To begin the construction of $\cA_n$, we let $\xi_0= {\bf j}$, and
consider the asymptotic lift $\gamma(t)$ of $\alfa(t)\in\cC $, given by
$$\pi(\gamma(t))=(\alfa(t),\alfa'(t)/||\alfa'(t)||),$$ where $\pi$ is as
in \eqref{doco}. We infer from this expression that the mapping
from $\cC$ into $C^{\8}([0,1])$ given by
\begin{equation}\label{exp1}
\alfa(t)\in\cC\mapsto
u_{\alfa}(t)=\int_0^t ||\gamma'(w)|| dw \in C^{\8}([0,1])
 \end{equation}
is continuous.

Next, note that for any fixed $\alfa(t)\in \cC$, the function
$u_{\alfa}(t)$ is non-negative and strictly increasing. Therefore we can define
the function $\hat{t}_{\alfa} (t):[0,1]\flecha [0,1]$ given by
 \begin{equation}\label{hatt}
\hat{t}_{\alfa}(t)=u_{\alfa}^{-1}\left(\frac{u_{\alfa} (t)}{n}\right),
 \end{equation}
which is smooth and strictly increasing. Hence we may consider the mapping
from
$\cC$ to
$C^{\8}([0,1])$ given by
 \begin{equation}\label{exp2}
\alfa(t)\in\cC\mapsto \hat{t}_{\alfa}(t) \in C^{\8}([0,1]).
 \end{equation}
We claim that this mapping is continuous. In fact, since the mapping
\eqref{exp1} is continuous, equation \eqref{hatt} shows that the
continuity of \eqref{exp2} would follow from the continuity of a mapping
$\cB$ that we describe next: let $\cS$ be the subset of
$C^{\8}([0,1])$ made up by all strictly increasing smooth functions
$f:[0,1]\flecha [0,+\8)$ such that $f(0)=0$. Then $\cB$ is the mapping from
$\cS$ into $\cS$ given by $$f(t)\in \cS\mapsto
f^{-1}\left(\frac{f(t)}{n}\right) \in \cS.$$ It is clear that
$\cB$ is well defined, because of the conditions imposed to $f(t)$.
Moreover, standard arguments show that $\cB$ is continuous, and therefore
\eqref{exp2} is also continuous, as we wished.

On the other hand, we can define another continuous mapping
from $\cC$ into $C^{\8}([0,1])$, namely the one given by
 \begin{equation}\label{exp3}
\alfa(t)\in \cC \mapsto k_{\alfa}(t) \in C^{\8}([0,1]),
 \end{equation}
where here $k_{\alfa}(t)$ stands for the geodesic curvature of $\alfa(t)$.
Next define, for a fixed element $\alfa(t)$ of $\cC$, the function
$\tilde{k}_{\alfa}(t):[0,1]\flecha \R$ given by $$\tilde{k}_{\alfa}(t)=
k_{\alfa}(\hat{t}_{\alfa} (t)).$$ Therefore we can define another mapping
from $\cC$ into $C^{\8}([0,1])$, this time given by
 \begin{equation}\label{exp4}
 \alfa (t)\in \cC \mapsto \tilde{k}_{\alfa}(t)\in C^{\8}([0,1])
 \end{equation}
Note that if $\alfa_1(t),\alfa_2(t)\in \cC$, the inequality
$$ \left|\tilde{k}_{\alfa_1}(t)-\tilde{k}_{\alfa_2}(t)\right|\leq
\left|k_{\alfa_1}(\hat{t}_{\alfa_1}(t)) - k_{\alfa_1}(\hat{t}_{\alfa_2}(t))\right| +
\left|k_{\alfa_1}(\hat{t}_{\alfa_2}(t))
-k_{\alfa_2}(\hat{t}_{\alfa_2}(t))\right|$$ holds.
As a consequence of this expression and the continuity of the maps \eqref{exp2} and
\eqref{exp3}, we obtain that \eqref{exp4} is again continuous.

Given $\alfa(t)\in\cC$, we plan to construct a curve
$\tilde{\alfa}(t):[0,1]\flecha\S^2$ verifying:

1) $u_{\alfa}(t)$ is the arclength parameter of the asymptotic lift in
$\S^3$ of
$\tilde{\alfa}(t)$, and

2) the geodesic curvature of $\tilde{\alfa}(t)$ is given by
$\tilde{k}_{\alfa}(t)$.

If $\tilde{s}_{\alfa}$ denotes the arclength parameter of a curve
$\tilde{\alfa}(t)$ in $\S^2$, then equation \eqref{relpar} ensures that
$\tilde{\alfa}(t)$ satisfies condition 1) if and only if
$\tilde{s}_{\alfa}(t)$ is given by
\begin{equation}\label{relpar2}
\frac{d\tilde{s}_{\alfa}}{dt}=\frac{d\tilde{s}_{\alfa}}{du_{\alfa}}
\frac{du_{\alfa}}{dt}=
\frac{2||\gamma'(t)||}{\sqrt{1+\tilde{k}_{\alfa}(t)^2}}.
\end{equation}
Note that the mapping
 \begin{equation}
 \alfa(t)\in \cC \mapsto \left(\frac{d\tilde{s}_{\alfa}}{dt},
 \tilde{k}_{\alfa}(t)\right)\in C^{\8}([0,1])\times C^{\8}([0,1]),
 \end{equation}
where $\tilde{s}_{\alfa}(t)$ is defined by \eqref{relpar2}, is trivially
continuous, and that so is also the process assigning to
$(d\tilde{s}_{\alfa}/dt,\tilde{k}_{\alfa}(t))$ the only (up to congruence)
curve $\tilde{\alfa}(t)$ in $\S^2$ with arclength parameter
$\tilde{s}_{\alfa}(t)$ and geodesic curvature $\tilde{k}_{\alfa}(t)$. This
follows from standard results on the regularity of the solutions to
ordinary differential systems with respect to initial conditions and
parameters.

Therefore, the mapping
 \begin{equation}\label{atil}
 \alfa(t)\in \cC \mapsto \tilde{\alfa}(t):[0,1]\flecha \S^2,
 \end{equation}
obtained via this process is well defined up to congruences in $\S^2$ and
continuous. Moreover, $\tilde{\alfa}(t)$ satisfies conditions 1) and 2).

Let us also note that if $\alfa(t)\in \cC$ and $k_1,k_2$ are smooth functions
satisfying that $k_1(u_{\alfa}(t))=\tilde{k}_{\alfa}(t)$, and
$k_2(u_{\alfa}(t))=k_{\alfa}(t)$, then the identity
 \begin{equation}
 k_1(u_{\alfa}(t))= k_2\left(\frac{u_{\alfa}(t)}{n}\right)
 \end{equation}
holds. This relation ensures that the angle $\omega(u)$ of the Hopf torus
$h^{-1}(\alfa(t))$, and the angle $\tilde{\omega}(u)$ of the Hopf cylinder
$h^{-1}(\tilde{\alfa}(t))$ are related by $\tilde{\omega}(u)=\omega(u/n)$.

To sum up, we have shown up to now that the process assigning to each
$\alfa(t)\in\cC$ the curve $\tilde{\alfa}(t):[0,1]\flecha\S^2$ such that,
with the above notations, $\tilde{\omega}(u)=\omega(u/n)$ holds, is well
defined up to congruences in $\S^2$, and continuous.

Next extend the parametrization $\alfa(t):[0,1]\flecha\S^2$ to a periodic
function $\alfa(t):\R\flecha\S^2$, so that $\alfa(t)$ (and hence
$\gamma(t)$) is $1$-periodic.

With this, observe that, since $k_{\alfa}(t)$ is $1$-periodic and
$\hat{t}_{\alfa} (t+n)= 1 +\hat{t}_{\alfa} (t)$,
the function $\tilde{k}_{\alfa}(t)$ is $n$-periodic. In this way
\eqref{relpar2} tells that $d\tilde{s}_{\alfa}/dt$ is also $n$-periodic
with respect to $t$, from where we obtain that
$\tilde{s}_{\alfa}(t)=\varrho (t) +c_0 t$, where $c_0\in \R$ and
$\varrho$ is $n$-periodic. But now, defining $\mu= c_0 n$, we get that
$\tilde{k}_{\alfa}(\tilde{s}_{\alfa}(t)
+\mu)=\tilde{k}_{\alfa}(\tilde{s}_{\alfa}(t))$, that is,
$\tilde{k}_{\alfa}$ is $\mu$-periodic with respect to the parameter
$\tilde{s}_{\alfa}$.

In order to finish the construction of the mapping $\cA_n$ from $\cC$ into
$\R$ mentioned at the beginning of the proof, let us consider for
$\tilde{\alfa}(\tilde{s})$ the only rigid motion $A_{\alfa}$ in $\S^2$
such that $A_{\alfa}(\tilde{\alfa}(0))=\tilde{\alfa}(\mu)$,
$A_{\alfa}(\tilde{\alfa}'(0))=\tilde{\alfa}'(\mu)$, and
$A_{\alfa}(\tilde{\alfa}(0)\times \tilde{\alfa}'(0))=
\tilde{\alfa}(\mu)\times \tilde{\alfa}'(\mu)$. It then holds that
$A_{\alfa}\in {\rm SO}(3)$, and moreover that $\tilde{\alfa}(\tilde{s}+
\mu)=A_{\alfa}(\tilde{\alfa}(\tilde{s}))$, since both curves have the same
geodesic curvature $\tilde{k}_{\alfa}(\tilde{s})$, and the same initial
conditions. Let $\theta_{\alfa}$ denote the angle of the rotation
$A_{\alfa}$, then it follows from the above comments that
$\theta_{\alfa}/\pi \in\Q $ if and only if $A_{\alfa}^q={\rm Id}$ for some
$q\in\N$, if and only if $\tilde{\alfa}(\tilde{s}+
q\mu)=\tilde{\alfa}(\tilde{s})$ for some $q\in \N$. In particular, if
$\tilde{k}_{\alfa}(\tilde{s})$ is not constant, we obtain that
$\tilde{\alfa}$ is closed if and only if $\theta_{\alfa}/\pi$ is rational.
Also note that the map
\begin{equation}\label{angu}
\tilde{\alfa}(t) \mapsto \theta_{\alfa}\in \R
\end{equation}
is continuous, and is invariant under congruences in $\S^2$.
Thus, we can define the continuous mapping
\begin{equation}\label{ar}
\cA_n:\cC\flecha \R
\end{equation}
obtained as the composition of \eqref{atil} and \eqref{angu}. This mapping
has been constructed so that it satisfies property {\bf (P)}, as desired.

It is known that there exist two equivalence classes of regular homotopy
on the set of smooth regular closed curves in $\S^2$. Those correspond to
the class of the circle, which will be denoted by $C_0$, and that of the
figure eight. One can check that if $\gamma$ is a circle in
$\S^2$, then $\cA_n(\gamma)=0\in \Q$. Actually, circles are the only
elements of
$\cC$ whose image under $\cA_n$ can be explicitly calculated. Since the mapping
$\cA_n:\cC\flecha \R$ is continuous, the set $J:=\cA_n (C_0)\subseteq\R$
is necessarily an interval
$J$, unless $\cA_n$ is constant on $C_0$, in which case $J=\{0\}$. In any
case, there exists an infinite number of curves $c(t)\in C_0$ that are not circles,
and such that $\cA_n(c(t))\in \Q$. Note that they can be chosen arbitrarily close
to any given circle. Finally, the Hopf torus of any of these
curves satisfies the conditions asserted in the Theorem, and we are done.
\end{proof2}
\begin{remark}
Since for every $n\in \N$, $n>1$ there exist infinitely many Hopf tori
$h^{-1}(c)$ such that $\cA_n (c)\in\Q$, and for every such $h^{-1}(c)$ we
have created a new flat tori in $\R^4$ via Theorem {\rm \ref{Boris}}, our
description actually produces a general family of new flat tori in $\R^4$,
and not just isolated examples.
\end{remark}
\begin{remark}
The construction process described in Theorems {\rm \ref{Boris}} and {\rm
\ref{truco}} can be seen as a method of unfolding a Hopf torus in
$\S^3$ of a particular kind, so that the result is a flat torus in $\R^4$
that does not lie in any affine
$3$-sphere.

To explain this interpretation, we will assume for simplicity that $n=2$.
So, we start with a Hopf torus $\Sigma_0=h^{-1}(c)$, with $c$ a closed
regular curve in $\S^2$, such that $\cA_2(c)\in\Q$, for the mapping
$\cA_2$ defined in \eqref{ar}. Thus $\cA_2(c)=l/k$, where $l,k\in \N$, and
this implies that the curve $\tilde{c}$ constructed from $c$ through
\eqref{atil} is closed. Moreover, $c(u)$ winds $k$ times whenever
$\tilde{c}(u)$ winds once, where here $u$ is the (common) arclength
parameter of their asymptotic lifts in $\S^3$. Next, let us consider a
trivial $2k$-folded covering of $\Sigma_0$, denoted by $\Sigma^*$,
obtained by tracing $c(u)$ $2k$-times, and taking then its Hopf torus in
$\S^3$. Thus $\Sigma^*$ is again a torus, and the Hopf torus
$\tilde{T}=h^{-1}(\tilde{c})$ can be parametrized
as a map from $\Sigma^*$ into $\S^3$ in the obvious way. Now, if we define
$\tilde{\alfa},\tilde{\beta}$ as coordinates of $\tilde{T}$ and its unit
normal in $\S^3$, and we parametrize them at \emph{double speed} in terms
of the asymptotic parameters $u,v$, i.e. we define
$(\alfa,\beta)(u,v)=(\tilde{\alfa},\tilde{\beta})(2u,2v)$, then
$\tilde{\alfa},\tilde{\beta}, \alfa,\beta$ are all well defined over
$\Sigma^*$. Once here we can use formula \eqref{rep} to obtain a map
$f:\Sigma^*\flecha\R^4$, which will be a flat torus in $\R^4$, possibly
with singular points, and not lying in any $3$-sphere. But finally, by
collapsing this map into $\Sigma_0$ we find that, closely enough to
$\Sigma_0$, $f$ does not have singular points.

In other words, we have covered $2k$-times the Hopf torus $\Sigma_0$, and
then perturbed this trivial covering in a special way, so that we assign
to each one of the $2k$-sheets of the covering a different perturbation.
After this deformation, the resulting surface is a regular flat torus in
$\R^4$, with flat normal bundle, regular Gauss map, and not lying in any
$3$-sphere, as we wished.
\end{remark}

A simpler version of the process we have just described can be used to
generate new complete flat cylinders in $\R^4$ with $R^{\perp}\equiv 0$
and regular Gauss map, as small perturbations of Hopf cylinders.

To begin with, let $N$ be a Hopf cylinder with angle $\omega(u)$. We shall
assume that $N$ is complete, has bounded mean curvature, and its
non-trivial family of asymptotic curves have bounded curvature in $\S^3$.

For $n\in \N$, $n>1$ we consider the Hopf cylinder $\tilde{N}$ with angle
$\tilde{\omega}(u)=\omega (u/n)$, as well as geometric solutions
$\tilde{\alfa}$, $\tilde{\beta}$ corresponding to $\tilde{\omega}$. Then
\eqref{ab} gives a solution of \eqref{sistema} for $\omega$, and we obtain
via \eqref{rep} a flat cylinder $f:u,v\text{-plane}\flecha\R^4$ with flat
normal bundle and regular Gauss map. It may have singular points, but it
does not lie in any $3$-sphere.

In order to get rid of singular points, we substitute $\alfa,\beta$ by
$\alfa_{\landa}=1+\landa \alfa$ and $\beta_{\landa}=\landa\beta$, and note
that
 \begin{equation}\label{converg}
\lim_{\landa\rightarrow 0} \left\{\left(A_{\landa}^2
-B_{\landa}^2\right)\sin\omega - 2A_{\landa}B_{\landa} \cos\omega
\right\}= \sin \omega
 \end{equation}
pointwise on the $u,v$-plane. Here $A_{\landa},B_{\landa}$ are the
functions appearing in \eqref{ambm}, but related to
$\alfa_{\landa},\beta_{\landa}$.

Besides, since $\tilde{\alfa},\tilde{\beta}$ are geometric solutions, the
flat surface in $\R^4$ that they define along with $\tilde{N}$ lies in a
$3$-sphere of $\R^4$. Thus, from \eqref{rep},
$\tilde{\alfa},\tilde{\beta},\tilde{\alfa}_u,\tilde{\beta}_u$ are bounded
on the $u,v$-plane. In addition, note that $\esiz
\tilde{N}_{uu},\tilde{N}_{uu}\esde = 1+\tilde{\omega}_u^2$. As
$\tilde{\alfa},\tilde{\beta}$ can be chosen to be coordinates of
$\tilde{N}$ and its polar map, we obtain
$\tilde{\alfa}_{uu}^2 +\tilde{\beta}_{uu}^2 \leq
2(1+\tilde{\omega}_u^2)$, and this ensures that
$\tilde{\alfa}_{uu}, \tilde{\beta}_{uu}$ are bounded, since we required
the asymptotic curves of $N$ to have bounded curvature.

Finally, note that
 $$\left\{\def\arraystretch{1.2} \begin{array}{lll}
 A_{\landa}(u,v) & = & 1 +\landa \left( \tilde{\alfa} +n^2 \tilde{\alfa}_{uu}
 +n^2\tilde{\omega}_u \tilde{\beta}_u\right)(nu,nv), \\
 B_{\landa}(u,v) & = & \landa \left( \tilde{\beta} +n^2 \tilde{\beta}_{uu}
 -n^2\tilde{\omega}_u \tilde{\alfa}_u\right)(nu,nv).
 \end{array}\right.$$
All of this shows that $A_{\landa}$ and $B_{\landa}$ are bounded. Hence,
the convergence in \eqref{converg} is actually uniform over the
$u,v$-plane. Since $N$ has bounded mean curvature, there is some $c_0$ such that
$0<c_0\leq \sin \omega$. Thus there is some $\landa_0>0$ such that, for all
$0<\landa<\landa_0$, $f_{\landa}$ preserves all the above mentioned properties of $f$, but
has no singular points.

Furthermore, the metric of $f_{\landa}$, given by \eqref{metrica} for
$A_{\landa},B_{\landa}$ and $\omega$, converges uniformly on the
$u,v$-plane to the metric of $N$, which satisfies $\esiz dN(X),dN(X)\esde
\geq d>0$ for some $d>0$, and for every unit tangent vector field
$X\in\X(N)$. This shows that, for
$\landa$ sufficiently small, $f_{\landa}:u,v\text{-plane}\flecha\R^4$ is
complete.

These are, to the best of our knowledge, the first examples of complete
flat cylinders in $\R^4$ with flat normal bundle and regular Gauss map
that do not lie in any $3$-sphere. Nevertheless, it is remarkable that all
of them are bounded in $\R^4$.

\def\refname{References}

\hspace{0.4cm}

\noindent J.A. Gálvez was partially supported by MCYT-FEDER, Grant no.
BFM2001-3318.

\hspace{0.2cm}

\noindent P. Mira was partially supported by MCYT, Grant no. BFM2001-2871
and CARM Programa Séneca, Grant no PI-3/00854/FS/01.


\begin{thebibliography}{999999}
\bibitem[Bia]{Bia} L. Bianchi, Sulle superficie a curvatura nulla
in geometria ellittica, {\it Ann. Mat. Pura Appl.}, {\bf 24} (1896),
93--129.

\bibitem[Bor]{Bor} A.A. Borisenko, Isometric immersions of space
forms in Riemannian and pseudo-Riemannian spaces of constant
curvature, {\it Russian Math. Surveys}, {\bf 56} (2001), 425--497

\bibitem[CaDa]{CaDa} M.P. do Carmo, M. Dajczer, Local isometric immersions
of $\R^2$ into $\R^4$, {\it J. Reine Angew. Math.},  {\bf 442}
(1993), 205--219.

\bibitem[CC]{CC} H. Cheng-Chung, A differential-geometric criterion for a
space curve to be closed, {\it Proc. Amer. Math. Soc.}, {\bf 83}
(1981), 357--361.

\bibitem[ChKa]{ChKa} C. Chicone, N.J. Kalton, Flat embeddings of the
Möbius strip in $\R^3$, {\it Comm. Appl. Nonlinear Anal.}, {\bf 9}
(2002), 31--50

\bibitem[DaTo]{DaTo} M. Dajczer, R. Tojeiro, On flat surfaces in space
forms, {\it Houston J. Math.}, {\bf 21} (1995), 319--338.

\bibitem[DaTo2]{DaTo2} M. Dajczer, R. Tojeiro, Submanifolds with nonparallel
first normal bundle, {\it Can. Math. Bull.}, {\bf 37} (1994),
330--337.

\bibitem[Eno]{Eno2} K. Enomoto, Global properties of the Gauss image of flat
surfaces in $R^4$, {\it Kodai Math. J.}, {\bf 10} (1987), 272--284.

\bibitem[GMM]{GMM} J.A. Gálvez, A. Martínez, F. Milán, Flat surfaces in the
hyperbolic $3$-space, {\it Math. Ann.}, {\bf 316} (2000),
419--435.

\bibitem[HKP]{HKP} T. Hasanis, D. Koutroufiotis, P. Pamfilos,
Surfaces of $E^4$ satisfying certain restrictions on their normal
bundle, {\it Trans. Amer. Math. Soc.}, {\bf 319} (1990), 329--347.

\bibitem[HoOs]{HoOs} D. Hoffman, R. Osserman, The geometry of the generalized
Gauss map, {\it Mem. Amer. Math. Soc.}, {\bf 28} (1980).

\bibitem[Kit1]{Kit1} Y. Kitagawa, Periodicity of the asymptotic curves on
flat tori in $S^3$, {\it J. Math. Soc. Japan}, {\bf 40} (1988),
457--476.

\bibitem[Kit2]{Kit2} Y. Kitagawa, Embedded flat tori in the unit
$3$-sphere, {\it J. Math. Soc. Japan}, {\bf 47} (1995), 275--296.

\bibitem[Pin]{Pin} U. Pinkall, Hopf tori in $\S^3$,  {\it Invent. Math.}, {\bf 81} (1985),
379--386.

\bibitem[Spi]{Spi} M. Spivak, {\it A comprehensive introduction to
differential geometry, Vol. IV}. Publish or Perish, Inc., Boston,
Mass., 1975.

\bibitem[Ten]{Ten} K. Tenenblat, {\it Transformation of manifolds
and applications to differential equations}. Pitman Monographs and
Surveys in Pure and Applied Mathematics. Longman, Harlow, 1998.

\bibitem[Wei1]{Wei1} J.L. Weiner, Flat tori in $\S^3$ and their Gauss
maps, {\it Proc. London Math. Soc.}, {\bf 62} (1991), 54--76.

\bibitem[Wei2]{Wei2} J.L. Weiner, Rigidity of Clifford tori, in
{\it Geometry and topology of submanifolds, VII}, World Sci.
Publishing, River Edge, NJ, 1995, pp. 274--277.

\bibitem[Wei3]{Wei3} J.L. Weiner, Isometric immersions of $\mathbf{E}^2$ into
$\mathbf{E}^4$, in {\it Geometry of submanifolds and related
topics} (Kyoto, 2001). Surikaisekikenkyusho Kokyuroku (2001), pp.
136--143.

\bibitem[Wns]{Wns} T. Weinstein, {\it An introduction to Lorentz
surfaces}.  Walter de Gruyter, Berlin, New York, 1996.

\bibitem[Wun]{Wun} W. Wunderlich, Über ein abwickelbares Möbiusband, {\it
Monatsh. Math}, {\bf 66} (1962), 276--289.

\bibitem[Yau]{Yau} S.T. Yau, Submanifolds with constant mean curvature II,
{\it Amer. J. Math.} {\bf 96} (1975), 76--100.
\end{thebibliography}
\end{document}